\documentclass{amsart}
\usepackage[T2A,T1]{fontenc}
\usepackage[utf8]{inputenc}
\usepackage[margin=1in]{geometry}

\newtheorem{theorem}{Theorem}[section]

\newtheorem{questions}{Open questions}[section]
\newtheorem{question}{Open question}[section]

\newtheorem{definition}{Definition}[section]
\newtheorem{conjecture}{Conjecture}[section]
\usepackage{tikz}
\usepackage{tikz-cd}
\DeclareMathOperator{\trace}{tr}
\newcommand{\A}{\mathcal{A}}
\newcommand{\matSP}[1]{\textsc{\MakeLowercase{#1}}}
\newcommand{\matSPgreek}[1]{#1}
\DeclareMathOperator{\Mat}{Mat}
\newcommand{\SP}{\mathrm{Sp}}
\newcommand{\SPA}{\SP_2(\mathcal{A},\sigma)}
\DeclareMathOperator{\End}{End}
\newcommand{\st}{\mid}
\newcommand{\sigmaT}[1]{\sigma(#1)}
\newcommand{\SL}{\mathrm{SL}}
\newcommand{\R}{\mathbb{R}}
\newcommand{\Adj}{\mathrm{Adj}}

\newcommand{\CC}{\mathbb{C}}
\newcommand{\HH}{\mathbb{H}}
\newcommand{\RR}{\mathbb{R}}

\newcommand{\ZZ}{\mathbb{Z}}

\newcommand{\KK}{\mathbb{K}}
\newcommand{\OO}{\mathbb{O}}

\newcommand{\SO}{\operatorname{SO}}

\newcommand{\Sp}{\operatorname{Sp}}

\newcommand{\Sym}{\operatorname{Sym}}
\newcommand{\Herm}{\operatorname{Herm}}

\newcommand{\lietype}[1]{\mathrm{#1}}

\renewcommand{\u}{\mathfrak{u}}

\newcommand{\CoordRing}{\mathcal{O}}

\tikzstyle{base}=[circle,draw,fill=black,inner sep=0pt, minimum width = 4pt]
\tikzstyle{frozen}=[rectangle,draw,fill=blue,inner sep=0pt, minimum width = 4pt,minimum height=4pt]
\tikzstyle{mutableBig}=[circle,draw=black,fill=black,inner sep=0pt,minimum size=10]
\tikzstyle{frozenBig}=[circle,draw=blue,fill=blue,inner sep=0pt,minimum size=10]
\tikzstyle{mutable}=[circle,draw=black,fill=black,inner sep=0pt,minimum size=6]
\tikzstyle{affine}=[circle,draw,fill=red,inner sep=0pt, minimum width = 4pt]
\tikzstyle{affine2}=[circle, draw, fill=red,inner sep=0pt, minimum width=8pt]
\tikzstyle{affine3}=[circle, draw, fill=red,inner sep=0pt, minimum width=12pt]
\tikzstyle{affine4}=[circle, draw, fill=red,inner sep=0pt, minimum width=16pt]
\tikzstyle{invis}=[circle,inner sep=0pt, minimum width=4pt]
\tikzstyle{fat2}=[circle,draw,fill=black,inner sep=0pt, minimum width=8pt]
\tikzstyle{fat3}=[circle,draw,fill=black,inner sep=0pt, minimum width=12pt]
\tikzstyle{fat4}=[circle,draw,fill=black,inner sep=0pt, minimum width=16pt]

\tikzstyle{vertex}=[circle,draw=black,fill=black,inner sep=0pt,minimum size=0]
\tikzstyle{line}=[thick]

\numberwithin{equation}{section}

\numberwithin{equation}{section}

\begin{document}

%------
% Insert the title of your paper and (if necessary)
% a short title for the running head.
%------

\title{Positivity and Non-commutativity}
\author{Anna Wienhard}

%------
% Insert your abstract.
%------
\begin{abstract}
In this article we revisit a new notion of positivity in real semisimple Lie groups that at the same time generalizes total positivity in split real Lie groups as well as positive Lie semigroups in Hermitian Lie groups of tube type. We shortly discuss the relationship with higher rank Teichm\"uller spaces, and then focus on describing different aspects of positivity as well as open questions. In the second part we describe a non-commutative perspective on Hermitian Lie groups of tube type that is suggested by positivity and leads to interesting applications, such as non-commutative generalizations of Markov numbers. 
\end{abstract}

\maketitle

%------
% INSERT THE BODY OF THE PAPER HERE (except
% acknowledgments, funding info and bibliography)
%------
\section{Introduction}
In my plenary lecture ``A beautiful world beyond hyperbolic geometry: Anosov representations and higher Teichmüller spaces'' I described aspects of hyperbolic geometry through the lens of groups acting by isometries on hyperbolic spaces, and gave a selective overview of recent developments. I focused particularly on Anosov representations of higher rank Teichm\"uller spaces. Anosov representations generalize many aspects of convex cocompact actions on hyperbolic spaces to more general discrete subgroups of Lie groups of higher rank. This leads to interesting group actions on higher rank symmetric spaces of non-compact type, such as the space of positive definite symmetric matrices or the Siegel upper half space. 
Special families of Anosov representations are given by higher rank Teichm\"uller spaces. Higher rank Teichm\"uller spaces are connected components of the representation variety of the fundamental group $\pi_1(S)$ of a closed surface into a semisimple Lie group, which consist entirely of discrete and faithful representations. For this to make sense, recall that classical Teichm\"uller space can be realized as the space of marked hyperbolic structures on $S$. The holonomy representations of marked hyperbolic structures form connected components of the representation vareity of $\pi_1(S)$ into the isometry group of the hyperbolic plane $\mathrm{PSL}(2,\mathbb{R})$. For  $\mathrm{PSL}(2,\mathbb{R})$ Anosov representations of 
 $\pi_1(S)$ are precisely the holonomy representations of hyperbolic structures on $S$. 
For other Lie groups, the existence of such higher rank Teichm\"uller spaces is in fact quite surprising, and in fact they only exist for special families of Lie groups. We refer the reader who is interested in learning more on Anosov representations to Canary's informal lecture notes on Anosov representations \cite{Canary_Anosov} and to Kassel's article in the Proceedings of the ICM \cite{Kassel_ICM}, and those who want to get a view into higher rank Teichm\"uller space to the author's article in the Proceedings of the ICM \cite{Wienhard_ICM}.

Higher rank Teichm\"uller spaces are intimately linked to positivity in Lie groups. In the case of split real Lie groups the connection to Lusztig's total positivity has been made by Fock and Goncharov \cite{FockGoncharov}. For Hermitian Lie groups of tube type, the connection to positivity has been made  in work of Burger, Iozzi, Wienhard \cite{BIW_annals}, characterisations in terms of causal structures and bi-invariant orders appeared in the work of Ben Simon, Burger, Hartnick, Iozzi and Wienhard \cite{BBHIW_order}. 

In recent years Guichard and Wienhard \cite{GW_pos1, GW_pos2} introduced a new notion of positivity in real semisimple Lie groups, which 
at the same time generalizes Lusztig's total positivity (defined for split real Lie groups), and positive Lie semigroups in Hermitian Lie groups of tube type. There are four families  of simple real Lie groups admitting a positive structure. Besides split real Lie groups and Hermitian Lie groups of tube type, one family is formed by indefinite orthogonal groups $\mathrm{SO}(p+1, p+k)$, where $p,k>1$. The other family consists of exceptional groups. 

Positivity gives rise to a notion of positivity of configurations of flags in appropriate flag varieties associated to $G$. This can be used to define positive representations of surface groups $\pi_1(S)$ into $G$ whenever $G$ admits a positive structures. It has been conjectured by Guichard, Labourie  and Wienhard \cite[Conjecture~19]{Wienhard_ICM} that the set of positive representations always form higher rank Teichmüller spaces. This has been proven recently by Beyrer, Guichard, Labourie, Pozzetti, and Wienhard \cite{GLW, BP, BGLPW}. Even though there are still several open questions on positive representations, this establishes that positivity provides the right underlying structure to explain the existence of higher rank Teichm\"uller spaces. 

But positivity is much more. It opens a new perspective on the families of simple real Lie groups admitting a positive structure. In this perspective non-commutative rings with involution play a crucial role. 
In this article we want to give a glimpse of that new perspective, and describe it in more details for Hermitian Lie groups of tube type. 
Through the lens of positivity, Hermitian Lie groups of tube type of arbitrary rank look like Lie groups of type $A_1$, but defined over a non-commutative ring with involution.  
This opens up a fresh view, which leads to interesting new constructions, including new quantum groups, generalizations of Markov numbers, non-commutative trace identities, and many other things which are yet to be explored. 

In the first part of this article we review important aspects of positivity in Lie groups and state several open questions. In the second part we illustrate the non-commutative perspective on Hermitian Lie groups of tube type and describe some of its consequences. 

%
%
% which is linked to non-commutativity in a crucial way. It gives rise to non-commutative generalizations of Markov triples, a fresh look onto the symplectic group and the Siegel upper half space, new non-commutative cluster algebras and motivates several conjectures.
%
%The goal of this article is to explore the new non-commutative perspectives, which positivity in Lie groups open. 

\section{Positivity in Lie groups} 
In this section we recall the notion of positivity in a real semisimple Lie group. It is a generalization of the notion of total positivity. 

The theory of totally positive matrices 
arose in the beginning of the 20th century through work of Schoenberg \cite{Schoenberg} and Gantmacher and Krein \cite{Gantmacher_Krein}. In the 1990's the theory has been generalized widely by Lusztig \cite{Lusztig}
who introduced total positivity in general split real semisimple Lie groups. Since then total positivity plays an important role in representation theory, cluster algebras, and has many relations to other areas in mathematics as well as to theoretical physics. 
The concept of positivity has been generalized by Guichard and Wienhard \cite{GW_pos1, GW_pos2} to more general real semisimple Lie groups which are not necessarily split. This notion of positivity is defined with respect to a subset 
$\Theta$ of the set of simple roots $\Delta$. The classification of semisimple Lie groups that admit a positive structure can be reduced to simple Lie groups. There are four families of simple  Lie groups that admit a positive structure. One family is the family of split real Lie groups with total positivity, which is positivity with respect to the full set $\Delta$. The three other families carry a positive structure with respect to a proper subset $\Theta$ of $\Delta$. 
As we will see below, when $\Theta$ is a proper subset of $\Delta$, some non-commutativity enters the picture. 
%Positive structures in Lie groups 
%
%- positivity defined in terms of the Lie algebra
%-positivity via positive triples in flag varieties. 

\subsection{The four families} 
Let us shortly recall the notion of positivity. For more details we refer the reader to \cite{GW_pos1, GW_pos2}. 

Consider a simple real Lie group $G$ which is not
necessarily split. A subset $\Theta\subset \Delta$ of the set of
simple roots defines a standard parabolic subgroup $P_\Theta$ in
$G$. (We choose the convention so that $P_\Delta$ is the minimal parabolic subgroup of~$G$.) The group $P_\Theta$ is % again
a semidirect product of a reductive Lie group $L_\Theta$ and its unipotent radical $U_\Theta$. 
The Lie algebra $\mathfrak{u}_\Theta$ of $U_\Theta$ carries a natural action
by $L_\Theta$ and can hence be decomposed into its $L_\Theta$-irreducible pieces. For every simple root $\alpha \in \Theta$ there is an irreducible piece $\mathfrak{u}_\alpha$. When $G$~is split and $\Delta = \Theta$ these would be precisely the root spaces~$\mathfrak{g}_\alpha$. In that case they are one-dimensional. 
When $\Theta$ is a proper subset of $\Delta$, then the $\mathfrak{u}_\alpha$ are in general not root spaces and not necessarily one-dimensional. 

The Lie group $G$ is said to admit a positive structure with respect to $\Theta$ if it satisfies the
following property:

\medskip
For every $\alpha \in \Theta$ there exists a non-empty closed acute convex
cone $c_\alpha \subset \mathfrak{u}_\alpha$ that is invariant under the action of $L_\Theta^\circ$ on $\mathfrak{u}_\alpha$, where $L_\Theta^\circ$ is the connected component of the identity. 
\medskip 

When $G$ is split, $\mathfrak{u}_\alpha =  \mathfrak{g}_\alpha$ can be identified with $\mathbb{R}$. Then the cone $c_\alpha$ can be identified with $\mathbb{R}^{\geq 0} \subset \mathbb{R}$. 

\begin{theorem}\cite[Theorem~3.4]{GW_pos1}
  A simple Lie group $G$ admits a positive structure with respect to $\Theta$ if and only if the pair $(G,\Theta)$ belongs to the following list: 
   \begin{enumerate}
  \item $G$ is a split real Lie group, and $\Theta = \Delta$;
  \item $G$ is Hermitian of tube type and of real rank~$r$ and $\Theta = \{ \alpha_r\}$, where $\alpha_r$ is the long simple restricted root;
  \item $G$ is locally isomorphic to $\mathrm{SO}(p+1,p+k)$, $p>1$, $k>1$ and
    $\Theta = \{ \alpha_1, \dots, \alpha_{p}\}$, where $\alpha_1, \dots,
    \alpha_{p}$ are the long simple restricted roots; %%% OG for "coherence"
                                %%% try to write p+1,p+k everywhere
  \item $G$ is the real form of $F_4$, $E_6$, $E_7$, or of $E_8$ whose reduced root system
    is of type~$F_4$, and $\Theta = \{ \alpha_1, \alpha_2\}$, where $\alpha_1,
    \alpha_{2}$ are the long simple restricted roots.
  \end{enumerate}
\end{theorem}

The arguments that give the classification provide important additional information. 

When $\Theta$ is a proper subset of $\Delta$, there exists a unique $\alpha_\Theta \in \Theta$, for which 
$\mathfrak{u}_{\alpha_\Theta}$ is not equal to the root space $\mathfrak{g}_{\alpha_\Theta}$ and is a higher dimensional vector space. 

For all other elements $\alpha$ of $\Theta \backslash\{\alpha_\Theta\}$, we have that 
$\mathfrak{u}_\alpha$ is equal to $\mathfrak{g}_\alpha$, and this root space has to be of dimension one. Hence, in this case, if we identify $\mathfrak{u}_\alpha$ with $\mathbb{R}$, the cone  $c_\alpha$ can be identified again with $\mathbb{R}^{\geq 0} \subset \mathbb{R}$.

\subsection{A non-commutative view on the special piece $\mathfrak{u}_{\alpha_\Theta}$}
For the special root $\alpha_\Theta$, the vector space $\mathfrak{u}_{\alpha_\Theta}$ carries an irreducible action of the reductive group $L_\Theta^{\circ}$, which preserves the convex cone $c_{\alpha_\Theta} \in \mathfrak{u}_{\alpha_\Theta}$ and acts transitively on its interior $\mathring{c}_{\alpha_\Theta}$. 
Such symmetric cones have been classified, they are intimately related to Hermitian Lie groups of tube type, see \cite{Faraut_Koranyi} for more details.

Thus let us look first at the case when $G$~is Hermitian of tube type. In this case there is just one special piece $\mathfrak{u}_{\alpha_\Theta}$.
$G$~is (up to
isogeny) $\mathrm{Sp}(2n, \R)$, $\mathrm{SU}(n,n)$, $\mathrm{SO}^*(4n)$, $\mathrm{SO}(2,q)$, or the exceptional group
$\mathrm{E}_{7(-25)}$. The $3$~first families of classical groups can be uniformly described: 
Let~$\KK$ be either~$\RR$, $\CC$, or~$\HH$, the space $V=\KK^{2n}$ is a right
vector space and is endowed with an anti-Hermitian form defined
by\index{$\omega_\KK$ the standard anti-Hermitian form on $\KK^{2n}$ ($\KK=\RR$,
  $\CC$, or $\HH$)}
\[
  \omega_\KK(v,v')= \sum_{i=1}^{n} \bar{x}_i y_{i}^{\prime} - \bar{y}_i x_{i}^{\prime}
\]
with $v=(x_1, \dots, x_n, y_1, \dots, y_n)$ and $v'=(x'_1, \dots, x'_n, y'_1, \dots, 
y'_n)$. The subgroup~$G$ in $\mathrm{GL}_\KK(V)$ of automorphisms of~$V$
stabilizing~$\omega_\KK$ is Hermitian of tube type; it is $\mathrm{Sp}(2n,\R)$,
$\mathrm{SU}(n,n)$, or $\mathrm{SO}^*(4n)$ following that $\KK$~is~$\RR$, $\CC$, or~$\HH$
respectively. 
The Levi factor~$L_\Theta$ in this case consists of the elements
\[ \ell_A =
  \begin{pmatrix}
    A & 0 \\ 0 & {}^t \overline{A}^{-1}
  \end{pmatrix}, \ A\in \mathrm{GL}_n(\KK ),
\]
and the unipotent algebra $\mathfrak{u}_\Theta$ is equal to
$\mathfrak{u}_{\alpha_\Theta}$. It consists of the matrices
\[
u_M  =
  \begin{pmatrix}
    0 & M \\ 0 & 0
  \end{pmatrix}, \ M\in M_n(\KK ) \text{ such that } {}^t\overline{M}=M;
\]
thus $\mathfrak{u}_{\alpha_\Theta}$~identifies with the space $\Herm(n,\KK)$ of Hermitian matrices over~$\KK$.
Since $ \ell_A u_M \ell_{A}^{-1} = u_{AM{}^t \overline{A}}$, we get that $\Herm^{\geq
  0}(n,\KK)$ is the sought for $L_\Theta$-invariant cone.

For
$\mathrm{SO}(2, k+1)$, let us consider the quadratic form\index{$q_{2,k+1}$ a
  quadratic form of signature $(2,k+1)$}
$q_{2,k+1}= -x_1 x_{k+3}+ x_2 x_{k+2} -\sum_{j=3}^{k+1} x_{j}^{2}$ on
$V=\R^{k+3}$. The group $G=\mathrm{SO}(q_{2,k+1})$ is also an Hermitian Lie group of tube type. 
%The Levi
%factor~$L_\Theta$ consists of the matrices
%\[
%  \ell_{\lambda, g} =
%  \begin{pmatrix}
%    \lambda \det g& 0 & 0 \\ 0 & g & 0 \\ 0 & 0 & \lambda^{-1}
%  \end{pmatrix},\ \lambda\in \RR^*, \ g\in \mathrm{O}(q_{1,k}),
%\]
%where $q_{1,k}$\index{$q_{1,k}$ a quadratic form of signature $(1,k)$} is given by restricting $q_{2,k+1}$ on the codimension~$2$
%space defined by $x_1=x_{k+3}=0$. 
The unipotent algebra is composed by the
matrices
\[
u_v =
\begin{pmatrix}
  0 & {}^t v J & 0 \\ 0 & 0 & v \\ 0 & 0 & 0
\end{pmatrix}, \ v\in \RR^{k+1},
\]
where $J$~is the matrix of the form~$q_{1,k}$. % , namely $J=
% \begin{pmatrix}
%   0 & 0 & 1 \\ 0 & \id_{k-1} & 0 \\ 1 & 0 & 0
% \end{pmatrix}
% $.
One has in this case $\ell_{\lambda,g} u_v \ell_{\lambda,g}^{-1} = u_{\lambda g v}$ and the cone
defined by $x_2>0$, $q_{1,k}>0$ is invariant by $L_{\Theta}^{\circ}$.

The exceptional Hermitian simple Lie group of tube type is $\mathrm{E}_{7(-25)}$ (or
EVII in Cartan's notation). In this case, the semisimple part of the
group~$L_\Theta$ is $\mathrm{E}IV =  \mathrm{E}_{6(-26)}$ and the action of~$L_\Theta$ on
$\mathfrak{u}_{\alpha_\Theta}$ is the action of the group of isomorphisms of $
\Herm(3,\OO)$ in it, $\Herm^{\geq 0}(3,\OO) $ is a natural invariant
cone (cf.\ \cite[p.~213]{Faraut_Koranyi}).

% \begin{table}[ht]
%   \centering
%  \begin{tabular}{lllll} 
% % \hline
% $G$ & $r$& $L_\Theta$& $\u_{\alpha_\Theta}$ & $c_{\alpha_\Theta}$  \\ 
%  \hline %\hline
%  $\Sp(2n,\R)$ &$n$& $\GL(n,\R)$ & $\Sym(n,\R)$ & $\Sym^{\geq 0}(n,\R)$  \\ 
% % \hline
%  $\SU(n,n)$ &$n$& $\GL(n,\C)$ & $\Herm(n,\CC)$ & $\Herm^{\geq 0}(n,\CC)$ \\
% % \hline
% $\SO^*(4n)$ &$n$& $\GL(n,\HH)$ & $\Herm(n,\HH)$ &$\Herm^{\geq 0}(n,\HH)$\\
% % \hline
% $\operatorname{EVII}=\lietype{E}_{7(-25)}$ &$3$& $\lietype{G}_2$& $ \Herm(3,\Oc)$ & $\Herm^{\geq 0}(3,\Oc) $\\
% % \hline
% $\SO(2,1+k)$ & $2$ & $\R^*\times \SO(1,k)$ & $(\R^{1,k},q_{1,k})$ & $\{ v \in \R^{1,k} \mid v_1\geq 0, q_{1,k}(v)\geq 0\}$\\ 
% % \hline
% \end{tabular}
%   \caption{The cones in the Hermitian cases; $r$~is the real rank of~$G$}
% \end{table}

\smallskip

In the general case, when $\Theta\neq \Delta$ and the
cardinality of~$\Theta$ is $>1$, 
the description of the cone $c_{\alpha_\Theta}$ reduces to the Hermitian
case. 
When $G = \SO(p+1, p+k)$, then % $\u_{\alpha_i}\cong \R$ for $1\leq i \leq p-1$,
(using standard numbering for the simple roots)
 $\mathfrak{u}_{\alpha_\Theta}\cong
(\R^{1,k},q_{1,k})$ with $c_{\alpha_\Theta} =  \{ v \in \R^{1,k} \mid
v_1\geq 0, q_{1,k}(v)\geq 0\}$ ($q_{1,k}$ denotes here the standard quadratic form
of signature~$(1,k)$ on~$\R^{1+k}$). 

When $G$ is one group in the exceptional family whose system of restricted roots is of type~$\lietype{F}_4$, then 
% $\mathfrak{u}_{\alpha_1} \cong \R$, and
$\mathfrak{u}_{\alpha_\Theta}$ and the cone
$c_{\alpha_2}$ are given in Table~\ref{tab:cones-F4}.

\begin{table}[ht]
  \centering
 \begin{tabular}{llll} 
% \hline
   $G$ & $\u_{\alpha_\Theta}$ &  $c_{\alpha_2}$  \\
   \hline
  %\hline
 $\operatorname{FI}=\operatorname{F}_{4(4)}$ & $\Sym(3,\R)$ & $\Sym^{\geq
                                                              0}(3,\R)$ {\rule{0pt}{2.6ex}}  \\ 
% \hline
 $\operatorname{EII}=\lietype{E}_{6(2)}$ & $\Herm(3,\CC)$ & $\Herm^{\geq 0}(3,\CC)$ \\
% \hline
$\operatorname{EVI}=\lietype{E}_{7(-5)}$ & $\Herm(3,\HH)$  &$\Herm^{\geq 0}(3,\HH)$\\
% \hline
$\operatorname{EIX}=\lietype{E}_{8(-24)}$ &$ \Herm(3,\OO)$  & $\Herm^{\geq 0}(3,\OO) $\\
% \hline
\end{tabular}
  \caption{The nontrivial cone for groups in the exceptional family}\label{tab:cones-F4}
\end{table}

There is a well known connection between symmetric cones and Euclidean Jordan algebras (see for example \cite{Faraut_Koranyi}). The Lie algebra $\mathfrak{u}_{\alpha_\Theta}$ carries the structure of a Euclidean Jordan algebra and the cone $c_{\alpha_\Theta}\in \mathfrak{u}_{\alpha_\Theta}$ that gives the positive structure is given as the cone of squares.

Here we want to stress a slightly different point of view. 
For the cases where $\mathfrak{u}_{\alpha_\Theta}$~identifies with the space $\Herm(n,\KK)$ of Hermitian matrices over~$\KK$, there is in fact an algebra $\A$ with an anti-involution  $\sigma$ such that $\mathfrak{u}_{\alpha_\Theta}$ identifies with $\A^{\sigma}$, the set of fix points of $\sigma$, and the cone is 
$c_{\alpha_\Theta} =  \{ a\sigma(a)\, |\, a \in \A\} \subset \A^{\sigma}$, and $\mathring{c}_{\alpha_\Theta} = \{ a\sigma(a)\, |\, a \in A \text{ invertible}\}$. 
In fact also the case when $\mathfrak{u}_\alpha  = \mathfrak{g}_\alpha$ fits in this setup. In this case $\A = \mathbb{R}$ and $\sigma$ is the identity map. 

The case of indefinite orthogonal groups, when $\mathfrak{u}_{\alpha_\Theta}\simeq \R^{(1,q+1)}$, a slightly modified picture holds where $\A$ is the Clifford group of a Clifford algebra. For the case when $\mathfrak{u}_{\alpha_\Theta}$~identifies with the space $\Herm(3,\OO)$, there is no ambient associative algebra, but many of the other features persists in that case as well, see for example \cite{GKNW2} for a detailed discussion. 

%
%In all of the above cases,
% there is a ring $\A$ and an (anti)-involution $\sigma$ such that $\u_{\alpha_\Theta} = \A^{\sigma}$, where $\A^{\sigma}$ denotes the set of fix points of $\sigma$. There are three cases 
% 
%\begin{enumerate}
%\item  The case $\A = \mathrm{Mat}_n(\mathbb{K})$ with $\mathbb{K} = \RR, \CC, \HH$, and if $n=3$ also $\mathbb{K} = \OO$. In the latter case this is not a ring, but can ignore this here. The involution $\sigma$ is the conjugate transpose ${}^t\overline{}$. 
%\item $\mathcal{R}$ is the even Clifford algebra associated to the vector space $V\simeq \R^{(1,q+1)}$, $\sigma(a) = e_1a^Te_1$ with $e_1$ the first element of the standard basis.
%\item In the cases when $\u_\alpha = \mathfrak{g}_\alpha$, $\A  = \mathbb{R}$ and $\sigma$ is the trivial involution. 
%\end{enumerate}
%
%In all of the above cases, the cone $c_{\alpha_\Theta} =  \{ a\sigma(a)\, |\, a \in \A\} \subset \A^{\sigma}$, and $\mathring{c_{\alpha_\Theta} = \{ a\sigma(a)\, |\, a \in \A invertible}$. 
%
%he Jordan product on $\u_{\alpha_\Theta}$ can be recovered f \cite{GKNW2}. 

%
% The square allows us to define the Jordan product $w\circ v$ on $\poscone{}$ by the formula $w\circ v = \frac{1}{2}((v+w)^2-v^2-w^2)$. We will not need this structure here. 

\subsection{The unipotent semigroup} 
An important role in positivity is played by the non-negative and positive semigroups of the unipotent group $U_\Theta$. 
The nonnegative semigroup $U_\Theta^{\geq 0}$ is the sub-semigroup of $U_\Theta$ generated by elements $\exp(v)$, with $v \in c_\alpha$ for some $\alpha \in \Theta$.
The positive semigroup $U_\Theta^{>0}$ is the interior of $U_\Theta^{\geq 0}$. 

 Similarly, one can define a nonnegative semigroup $U_\Theta^{\mathrm{opp},\geq 0}$ and and a positive semigroup $U_\Theta^{\mathrm{opp},>0}$ in $U_\Theta^{\mathrm{opp}}$. 
This gives rise to the nonnegative semigroup $G_\Theta^{\geq 0}$ in $G$, which is the semigroup generated by  $U_\Theta^{\geq 0}$,  $U_\Theta^{\mathrm{opp},\geq 0}$, and $L_\Theta^\circ$. 
The positive semigroup $G_\Theta^{>0}$ is the interior of $G_\Theta^{\geq0}$.

In fact, for the reader who is not happy with the definition of positivity, requiring the existence of invariant convex cones, the existence of such semigroups could alternatively be used as the definition of positivity. 
\begin{theorem}\cite{GW_pos2}
 Let~$G$ be a connected simple Lie group.  
 Suppose that there is $U^{+}\subset U_\Theta$ such that 
 $U^{+}$ is a closed $L_{\Theta}^{\circ}$-invariant semigroup of non-empty interior, which contains no nontrivial invertible element. 
  Then $G$ admits a positive structure with respect to $\Theta$ and the semigroup~$U^+$ contains the
  semi\-group~$U_{\Theta}^{\geq 0}$.
  
  If we further assume that the interior of $U^+$ is contained in the open Bruhat cell with respect to $P_\Theta^{\mathrm{opp}}$,
  then $U^+ = U_{\Theta}^{\geq 0}$.
\end{theorem}
 
 Defining the positive semigroup $U_\Theta^{>0}$ as the interior of $U_\Theta^{\geq 0}$ is easy to state but does not provide much insight in how elements in this semigroup look like. 

In Lusztig's work, the positive semigroup $U_\Delta^{>0}$ is in fact defined by giving explicit parametrizations, proving that the image of the parametrizations do not depend on the choices made. The thus defined set coincides with the interior of  $U_\Delta^{\geq 0}$ is proved afterwards. The parametrizations are associated to reduced expressions of the longest element in the Weyl group, and the change of coordinates between parametrizations associated to different reduced expressions is given by an explicit formula involving only positive rational functions.  

In the case of positive structures with respect to a proper subset $\Theta$ of $\Delta$ similar results hold. However the Weyl group has to be replaced by another Coxeter group, called the $\Theta$-Weyl group. Every reduced expression of the longest word in the $\Theta$-Weyl group gives rise to a parametrization of the positive semigroup $U_\Theta^{>0}$, and the change of coordinates between two different reparametrizations is determined by explicit formulas \cite{GW_pos1,GW_posneu}.

% In the general case of positivity with respect to a proper subset $\Theta$ of $\Delta$, parametrizations of the positive semigroup $U_\Theta^{>0}$ are governed by the $\Theta$-Weyl group. 
% 
\subsection{The $\Theta$-Weyl group} 
%A crucial role in understanding Lie groups with a positive structure with respect to $\Theta$ is played by the $\Theta$-Weyl group $W(\Theta)$. 
The $\Theta$-Weyl group is a subgroup of the Weyl group $W$ of $G$. It turns out to be the Weyl group of a different root system, and lies at the heart of regarding the Lie groups admitting a positive structure with respect to $\Theta$ through the lens of non-commutativity. 

To introduce the $\Theta$-Weyl group, let us recall that the Weyl group $W$ is generated by reflections (order two elements) $s_\alpha$, $\alpha \in \Delta$. These elements are characterized by their action on the collection of roots or root spaces. The element $s_\alpha$ sends the root space $\mathfrak{g}_\alpha$ to $\mathfrak{g}_{-\alpha}$ and permutes all the other root spaces associated to positive roots. 

The $\Theta$-Weyl group
$W(\Theta)$ is a subgroup of $W$ generated by elements $\sigma_\alpha$ for all
$\alpha \in \Theta$. For all $\alpha \in \Theta$ which are not connected to
$\Delta\setminus\Theta$, i.e. when $\mathfrak{u}_\alpha$ is equal to $\mathfrak{g}_\alpha$,  $\sigma_\alpha$ is just the reflection~$s_\alpha$. For the unique $\alpha_\Theta\in \Theta$ which  is connected toWeyl
$\Delta\setminus\Theta$, i.e. for which $\mathfrak{u}_{\alpha_\Theta}$ is a higher dimensional vector space, the element $\sigma_{\alpha_\Theta}$ is a specific word in the subgroup
generated by  the reflections associated to roots in $ \{\alpha_\Theta\} \cup \Delta\setminus\Theta$. The precise word is not of importance to us here, the reader can consult \cite{GW_pos1}.
The key point is that the action of the $\Theta$-Weyl group on the decompostion of the Lie algebra 
\[
\mathfrak{g} = \oplus_{\alpha \in <\Theta>} \mathfrak{u}_\alpha \oplus \mathfrak{l}_\Theta \oplus  \oplus_{\alpha\in <\Theta>} \mathfrak{u}_{-\alpha} 
\]
is completely analogous to the action of the classical Weyl group on the root decomposition 
\[
\mathfrak{g} = \oplus_{\alpha \in <\Delta>} \mathfrak{g}_\alpha \oplus \mathfrak{g}_0 \oplus  \oplus_{\alpha \in <\Delta>} \mathfrak{g}_{-\alpha}. 
\]

In other words, elements $\sigma_\alpha$ act on the irreducible pieces $\mathfrak{u}_\alpha$ in the following way. The element $\sigma_\alpha$ sends  $\mathfrak{u}_\alpha$ to  $\mathfrak{u}_{-\alpha}$ and permutes the $\mathfrak{u}_\beta$ for all the other $\beta$ in the positive span of $\Theta$. 

The $\Theta$-Weyl group is itself an abstract Weyl group, but of a different type than the Weyl group of $G$. 
Its rank is the cardinality of $\Theta$, and it acts on the pieces $\mathfrak{u}_\alpha$ in the same way as the the abstract Weyl group would act on the classical root spaces. 

Here is a list of the $\Theta$-Weyl groups. 
   \begin{enumerate}
  \item {\bf split case}: If $G$ is a split real form and $\Theta = \Delta$, then $W(\Theta) = W$
  \item {\bf $A_1$ case}: If $G$ is Hermitian of tube type and of real rank~$r$ and $\Theta = \{ \alpha_r\}$, then the Weyl group is of type $C_r$, but the $\Theta$-Weyl group is of type $A_1$. 
   \item {\bf $B_p$ case}:  If $G$ is locally isomorphic to $\mathrm{SO}(p+1,p+k)$, $p>1$, $k>1$ and
    $\Theta = \{ \alpha_1, \dots, \alpha_{p}\}$, then the Weyl group is of type $B_{p+1}$, but the $\Theta$-Weyl group is of type $B_p$.
  \item {\bf $G_2$ case}: If $G$ is the real form of $F_4$, $E_6$, $E_7$, or of $E_8$ whose reduced root system
    is of type~$F_4$, and $\Theta = \{ \alpha_1, \alpha_2\}$, then the Weyl group is of type $F_4$, but the $\Theta$-Weyl group is of type $G_2$.
      \end{enumerate}

The labels we assigned to the different cases reflect the fresh perspective we suggest to take when looking at the groups $G$ admitting a positive structure with respect to $\Theta$. 
We suggest to consider $G$ as Lie groups of type $W(\Theta)$, but defined over a (partially) non-commutative ring with an involution. 
We will illustrate this perspective below for the case of Hermitian Lie groups of tube type. In this case $\Theta = \{\alpha_\Theta\}$, so the $\Theta$-Weyl group is $\mathbb{Z}/2\mathbb{Z}$. Thus we want to consider Hermitian Lie groups of tube type as groups of type $A_1$, or more precisely as symplectic groups $\mathrm{Sp}_2$ over some non-commutative ring. 
This has been made precise in \cite{ABRRW, Rogozinnikov_20} and leads to interesting new phenomena. 

Before we do illustrate this non-commutative perspective, we describe how the $\Theta$-Weyl group appears in parametrizations of the positive semigroup, and state some open questions on parametrizations of the non-negative semigroups. 

\subsection{Parametrizing the positive semigroup}
We explain now how reduced expressions of the longest work in the $\Theta$-Weyl group give parametrizations of the positive semigroups $U_\Theta^{>0}$. 

Given any reduced expression  $\mathbf{\gamma} =\sigma_{\gamma_1} \cdots \sigma_{\gamma_n}$ of the
longest element $w_{\max}^{\Theta} \in W(\Theta)$ we define the map \[F_{\mathbf{\gamma}}\colon
\mathfrak{u}_{\gamma_1} \times \cdots \times \mathfrak{u}_{\gamma_n} \to U_\Theta\] as the
product of the exponential map on each factor, i.e. 
\[ F_{\mathbf{\gamma}} (X_1, \cdots, X_n)  = \exp(X_1) \cdots \exp(X_n),\]
where $X_i \in \mathfrak{u}_{\gamma_i}$.

Since $\gamma_i \in \Theta$ for all $i$, we can restrict the map $ F_{\mathbf{\gamma}} $ to the product of the cones  $c_{\gamma_i} \subset \mathfrak{u}_{\gamma_i}$ as well as to their interiors $\mathring{c}_{\gamma_i}$. 
The key properties can be summarized by the following 

\begin{theorem}\cite[Section~10]{GW_pos1}
\begin{enumerate}
\item The map $ F_{\mathbf{\gamma}} $ restricted to the product of open cones $\mathring{c}_{\gamma_1} \times \cdots \times \mathring{c}_{\gamma_n}$ is injective and $F_{\mathbf{\gamma}} (\mathring{c}_{\gamma_1} \times \cdots \times \mathring{c}_{\gamma_n}) = U_\Theta^{>0}$. Thus $F_{\mathbf{\gamma}} $ gives a parametrization of the positive semigroup. 
%\item $F({\mathring{c}_{\gamma_1} \times \cdots \times \mathring{c}_{\gamma_n}})$ is contained in the open Bruhat cell $\Omega_\Theta^{\mathrm{opp}}$ with respect to $P_\Theta^{\mathrm{opp}}$. In fact, it is a connected component of $\Omega_\Theta^{\mathrm{opp}} \cap U_\Theta$.
\item  The map $ F_{\mathbf{\gamma}} $ restricted to the product of closed cones ${c}_{\gamma_1} \times \cdots \times {c}_{\gamma_n}$ is not injective, but proper,  and $F_{\mathbf{\gamma}} ({c}_{\gamma_1} \times \cdots \times{c}_{\gamma_n}) = U_\Theta^{\geq 0}$.
\end{enumerate} 
\end{theorem}

Note that in the case of split real Lie groups all the cones $\mathring{c}_{\gamma_i} \subset c_{\gamma_i}\subset \mathfrak{u}_{\gamma_i}$ can be identified with $\mathbb{R}^{>0}\subset\mathbb{R}^{\geq 0}
 \subset \mathbb{R}$. Given two different reduced expression ${\mathbf{\gamma}}$ and ${\mathbf{\gamma'}}$ the change of coordinates are given by positive rational functions \cite{Lusztig}. Since any two different reduced expressions are related by a sequence of braid relations, it suffices to establish that for the braid relations in all split real Lie groups of rank 2. Berenstein and Zelevinsky \cite{BZ} wrote explicit expression for the change of coordinates. 
To derive these expression, given a braid relation, they set up a system of polynomial equations in the universal enveloping algebra, which they then solve explicitly. 
 
In the case when $\Theta$ is a proper subset of $\Delta$, several of the cones in the domain of $ F_{\mathbf{\gamma}} $ are cones in higher dimensional vector spaces. By the above theorem, the image of $F_{\mathbf{\gamma}}$ is still independent of the reduced expression $\mathbf{\gamma}$ of the longest word in the $\Theta$-Weyl group $W(\Theta)$. 
In \cite{GW_posneu} we give explicit expressions for the change of coordinates. These expressions are not positive rational functions, but surprisingly they still preserve the positivity defined by the cones. In fact, associated to braid relations in the $\Theta$-Weyl group we derive an analogous system of ``polynomial equations'' as Berenstein and Zelevensky, just that some of the variables are now not scalars, but take values in a vector space (or better a Euclidean Jordan algebra), and the powers of such an element are interpreted in an appropriate way as taking a product or evaluating a linear form. For details we refer the reader to \cite{GW_posneu}. 

What we would like to highlight here - where in the split real case,  the system of polynomial equations is associated to the type of the Weyl group, in the general case of positive structures with respect to $\Theta$ it is the same system of polynomial equations associated to the type of the $\Theta$-Weyl group, where the variable associated to the special root in $\Theta$ should be thought of as elements in a Euclidean Jordan algebra. 

Lusztig's foundational work on total positivity in Lie groups has spurred a lot of further research. Many of the results that have been established for total positivity in split real Lie groups, lead to natural questions for positivity with respect to $\Theta$ in more general semisimple real Lie groups. In the following subsections we mention some open questions, but literally almost any result you come across for total positivity, you can ask what the corresponding statement would be for positivity with respect to $\Theta$. 

\subsection{Parametrizations of the non-negative  semigroups}
The map $ F_{\mathbf{\gamma}} $ associated to a reduced expression of the longest word in the $\Theta$-Weyl group gives parametrization of the positive semigroup $U_\Theta^{> 0}$. This parametrization is equivariant with respect to the action of $L_\Theta^\circ$. 
 Since the map $ F_{\mathbf{\gamma}} $ is not injective when restricted to the product of closed cones ${c}_{\gamma_1} \times \cdots \times {c}_{\gamma_n}$, it does not give a parametrization of the nonnegative semigroup $U_\Theta^{> 0}$. 

To approach the problem of parametrizing the nonnegative semigroup it is instructive to look at the Bruhat decompositions of a Lie group $G$. 
Recall  that $G$~can be decomposed under the left-right multiplication
by $P_{\Delta}^{\mathrm{opp}}\times P_{\Delta}^{\mathrm{opp}}$, where
$P_{\Delta}^{\mathrm{opp}}$ is the minimal standard opposite parabolic subgroup.  These orbits are
 called the Bruhat cells and are indexed by elements of the Weyl group $W$.
One can also consider decompositions with respect to the action of $P_{\Delta}^{\mathrm{opp}}\times P_\Theta^{\mathrm{opp}}$ and $P_\Theta^{\mathrm{opp}}\times P_\Theta^{\mathrm{opp}}$, where $\Theta$~is a subset of~$\Delta$. For these the orbits are indexed by elements of  $W/W_{\Delta-\Theta}$, respectively $W_{\Delta-\Theta}\backslash W/W_{\Delta-\Theta}$, where $W_{\Delta-\Theta}$ is the Weyl group of $L_\Theta$. 

\begin{enumerate}
%\item \label{item:5:Bruhat_dec} For every~$\alpha$ in~$\Delta$,
%  \[P_{\Delta}^{\mathrm{opp}} s_\alpha P_{\Delta}^{\mathrm{opp}} s_\alpha P_{\Delta}^{\mathrm{opp}} =
%    P_{\Delta}^{\mathrm{opp}}  \sqcup P_{\Delta}^{\mathrm{opp}} s_\alpha P_{\Delta}^{\mathrm{opp}}.\]
%\item For every~$w$ in~$W$, the $P_{\Delta}^\mathrm{opp}\times
%  P_{\Theta}^{\mathrm{opp}}$-orbit
%  \[ P_{\Delta}^\mathrm{opp} w
%    P_{\Theta}^{\mathrm{opp}}\]
%  depends only on the class~$x$ of $w$ in $W/ W_{\Delta\setminus\Theta}$
%  and
%  will sometimes be denoted $P_{\Delta}^\mathrm{opp} x
%  P_{\Theta}^{\mathrm{opp}}$.\index{$P_{\Delta}^\mathrm{opp} x
%  P_{\Theta}^{\mathrm{opp}}$ the $P_{\Delta}^\mathrm{opp} \times
%  P_{\Theta}^{\mathrm{opp}}$-orbit associated with the class $x\in
%  W/W_{\Delta\smallsetminus \Theta}$} 
\item The group~$G$ is the disjoint union
  \[ G = \bigcup_{{w\in W}} P_{\Delta}^\mathrm{opp}w
  P_{\Delta}^{\mathrm{opp}}.\]
\item The group~$G$ is the disjoint union
  \[ G = \bigcup_{{x\in W/W_{\Delta\setminus\Theta}}} P_{\Delta}^\mathrm{opp}x
  P_{\Theta}^{\mathrm{opp}}.\]
\item The group~$G$ is the disjoint union
  \[ G =\!\!\!\!\!\!\!\!\!\bigcup_{a\in W_{\Delta\setminus\Theta}\backslash W
      /W_{\Delta\setminus\Theta}}\!\!\!\!\!\!\!\!\! P_{\Theta}^{\mathrm{opp}} a
  P_{\Theta}^{\mathrm{opp}}.\]
\end{enumerate}

In the case of total positivity, when $\Theta = \Delta$  all these decompositions agree. The intersection $U_{w}^{\geq 0}$ of the nonnegative semigroup $U^{\geq 0}$ with the Bruhat cell 
$P_{\Delta}^\mathrm{opp}w
  P_{\Delta}^{\mathrm{opp}}$ 
  can be parametrized in a similar way as the positive semigroup \cite{Lusztig}. Let ${\mathbf{\gamma} }= \gamma_1 \cdots \gamma_k$ be a reduced expression of an element $w \in W$, not necessarily the longest element . Then $F_{\mathbf{\gamma}}$ restricted to 
$({\mathring{c}_{\gamma_1} \times \cdots \times \mathring{c}_{\gamma_k}})$ gives a parametrization of $U_{w}^{\geq 0}$. Clearly $U^{\geq 0} = \bigcup_{{w\in W}} U_{w}^{\geq0}$ is a disjoint union. Thus this gives a parametrization of $U^{\geq 0}$, and in fact also a cell decomposition \cite{Rietsch}. 

In the case when $\Theta$ is a proper subset of $\Delta$ the task to parametrize the non-negative semigroup becomes more difficult. 
This becomes already apparent in the case of Hermitian Lie groups of tube type, which is the non-commutative $A_1$-type situation, where $\Theta = \{ \alpha\}$ consists of a single simple root $\alpha$. 
Here the situation is simple, there is just one space $\mathfrak{u}_{\alpha}$ and cones $\mathring{c}_{\alpha} \subset c_{\alpha}\subset \mathfrak{u}_{\alpha}$. The exponential map 
$\mathfrak{u}_{\alpha} \rightarrow U_{\{\alpha\}}$ is an isomorphism, hence the map $F: c_{\alpha} \rightarrow U_{\{\alpha\}}^{\geq 0}$ gives a parametrization. 

However, the structure of the intersection  of $U_{\{\alpha\}}^{\geq 0}$ with the different orbits in the above decompositions of $G$ is already very interesting. 
This steams from the fact that the boundary of the cone $\partial  c_{\alpha} = c_{\alpha}\ - \mathring{c}_{\alpha}$ is non trivial. In the split real case, or when $\mathfrak{u}_{\alpha}$ is of dimension one, this boundary is just $0$, but in general it has several strata, given by  $L_\Theta^\circ$ invariant subsets. These strata are precisely the orbits $U_{\{\alpha\}}^{\geq 0} \cap  P_{\Theta}^{\mathrm{opp}} a
  P_{\Theta}^{\mathrm{opp}}$ for $a\in W_{\Delta\setminus\Theta}\backslash W /W_{\Delta\setminus\Theta}$, see \cite{GW_pos1}. Note that this decomposition does not give you a cell decomposition as the individual strata have non-trivial topology. 
  
  The decompositions of $U_{\{\alpha\}}^{\geq 0}$ with respect to the other decompositions $G = \bigcup_{{x\in W/W_{\Delta\setminus\Theta}}} P_{\Delta}^\mathrm{opp}x
  P_{\Theta}^{\mathrm{opp}}$ and $G = \bigcup_{{w\in W}} P_{\Delta}^\mathrm{opp}w
  P_{\Delta}^{\mathrm{opp}}$, are more involved, as they are not invariant by $L_\Theta^\circ$. 
  The paper \cite{GW_pos1} contains further information with respect to the first decomposition, see also \cite{Xie, Xie_Wienhard} for further work in the case of the symplectic group. In particular for the symplectic group \cite{Xie} gives a cell decomposition of the non-negative semigroup $U_{\{\alpha\}}^{\geq 0}$. 
  
%   
%  In the special case of the symplectic group, the article \cite{Kaitao} gives a quite detailed description of the first decomposition, a description with respect to the  classical Bruhat decomposition $G = \bigcup_{{w\in W}} P_{\Delta}^\mathrm{opp}w P_{\Delta}^{\mathrm{opp}}$ will be describe in forthcoming work \cite{Kaitao_Wie}. 

For $\mathrm{SO}(p+1, p+k)$  where the $\Theta$-Weyl group is if type $B_p$, an explicit parametrization of the nonnegative semigroup $U_\Theta^{\geq0}$ has been given in  \cite[Section~11]{GW_pos1}, but even in this case the full structure of parametrizations of the nonnegative semigroup is not known. 
There are several open questions regarding the non-negative semigroups $U_\Theta^{\geq 0}$

\begin{questions}
\begin{enumerate}
\item Describe the strata $U_{\Theta, a}^{\geq 0}$, $U_{\Theta,x}^{\geq 0}$ and $U_{\Theta,w}^{\geq 0}$ of the nonnegative semigroup $U_\Theta^{\geq 0}$ with respect to the different decompositions of $G$ above. 
\item Provide parametrizations of $U_\Theta^{\geq 0}$ and of the strata  $U_{\Theta, a}^{\geq 0}$, $U_{\Theta,x}^{\geq 0}$ and $U_{\Theta,w}^{\geq 0}$. 
\item Find explicit cell decompositions of $U_\Theta^{\geq 0}$. 
\end{enumerate}
\end{questions}

\subsection{Positive parts of flag varieties}
Totally positive and nonnegative semigroups give rise to positive or non-negative parts in flag varieties. The positive Grassmannian has received a lot of interests, not only in the mathematical community, but also in theoretical physics, see Williams ICM talk \cite{Williams_ICM} for an overview.

If $G$ admits a positive structure with respect to $\Theta$, the non-negative semigroup $G_\Theta^{\geq 0}$ and the positive semigroup $G_\Theta^{>0}$ also give rise to non-negative parts and positive parts of flag varieties of a group G. These are introduced and studied in \cite{GW_posneu}. 

Let us give a short description here. 
Let $\mathcal{F} = G/P$ be a flag variety of $G$. 
The non-negative part of the flag variety  
$\mathcal{F}^{\geq0}$  with respect to $\Theta$ is 
\[{\mathcal F}^{\geq 0} := \{  G_\Theta^{\geq 0} P \} =\{ uPu^{-1} \,|\,u \in G_\Theta^{\geq0}\} ,\] where $P$ is a parabolic subgroup compatible with $P_\Theta$, i.e. such that $P \cap P_\Theta$ is a parabolic subgroup. 

The positive part $\mathcal{F}^{>0}$ of the flag variety ${\mathcal F}$ is 
\[
{\mathcal F}^{>0} := \{  G_\Theta^{>0} P \} = \{ uPu^{-1} \,|\,u \in G_\Theta^{>0}\},\] where $P$ is a parabolic subgroup compatible with $P_\Theta$. 
The non-negative part of the flag variety  $\mathcal{F}^{\geq0}$ is the closure of $\mathcal{F}^{>0}$ in $\mathcal{F}$. 

If $P= P_\Theta$ the positive part of $\mathcal{F}$ is what is called a diamond in \cite{GW_pos1, GLW}, and the non-negative part is the closure of a diamond. 

In the case of the symplectic group (with the $\Theta$-Weyl group being of type $A_1$), the positive and non-negative parts of the flag varieties have a nice geometric description, see  \cite{GW_posneu}.
The non-negative part of the Lagrangian Grassmannian is studied in detail in \cite{Xie}. In general, the structure of the non-negative and positive parts of flag varieties when $G$ is not necessarily split and the $\Theta$-Weyl group is of rank at least two, is not yet well understood.  

Here are some  open problems  
\begin{questions}
\begin{enumerate}
\item Provide explicit parametrizations of the positive and non-negative parts of flag varieties. 
\item Give cell decompositions of  the positive and non-negative parts of flag varieties. In the case of total positivity, for Grassmanians these are parametrized by positroid varieties. 
\item When $G$ is a split real form, which admits another positive structure with respect to a proper subset $\Theta$ in $\Delta$, describe precisely the relation between the positive and non-negative parts of the flag variety, with respect to $\Delta$ and with respect to $\Theta$. 
\item The amplituhedron, introduced by \cite{ArkaniHamed_Trnka}, is the image of the positive Grassmannian under a positive linear map. It plays a major role in understanding of scattering amplitudes in gauge theories. 
Introduce a generalization of the amplituhedron for positive structures with respect to $\Theta$. 
\end{enumerate}
\end{questions}

In the totally positive case, the positive parts of flag varieties can be described by cluster coordinates, so we expect that the non-commutative cluster structures associated to groups admitting a positive structure, combined with a good understanding of the situation when the $\Theta$-Weyl group if of type $A_1$ should be useful to answer the above questions. 

Thus now we turn to describe in more detail the non-commutative perspective on Lie groups admitting a positive structure with respect to $\Theta$, with a special focus on the case when the $\Theta$-Weyl group is of type $A_1$. For this we will first recall some aspects of the simple real Lie group of type $A_1$, namely, $\mathrm{SL}_2(\mathbb{R})$.

\section{Aspects of $\mathrm{SL}_2(\mathbb{R})$}\label{sec:SL2}

The Lie group $\mathrm{SL}_2(\mathbb{R})$ is the simplest real simple Lie group. It (and its Lie algebra) plays an important role in understanding the structure of any real simple Lie group. In this section we just highlight some aspects of 
$\mathrm{SL}_2(\mathbb{R})$ which will play a role later when we consider Lie groups of type $\mathrm{SL}_2$ over non-commutative rings. 

\subsection{$\mathrm{SL}_2(\mathbb{R})$ as $\mathrm{Sp}_2(\mathbb{R})$}
Usually we think of $\mathrm{SL}_2(\mathbb{R})$ as the subgroup of the general linear group $\mathrm{GL}_2(\mathbb{R})$ of determinant one. However, here we would like to think of $\mathrm{SL}_2(\mathbb{R})$ as 
$\mathrm{Sp}_2(\mathbb{R})$, the subgroup invertible endomorphism preserving a symplectic form $\omega$

The symplectic form $\omega$ is  a non-degenerate skew-symmetric bilinear form 
\[
\omega:\mathbb{R}^2 \times \mathbb{R}^2 \rightarrow \mathbb{R}, 
\]
defined by $\omega ((x_1, x_2), (y_1, y_2) )  = x_1 y_2 - x_2y_1.$
%\item As an element $\omega$ of the second exterior power $\Lambda^2 {\mathbb{R}^2}^*$, which is given by the tensor $\omega = {e_1}^* \otimes {e_2}^*  - {e_2}^*\otimes {e_1}^*$, where $e_1, e_2$ is the standard basis of $\mathbb{R}^2$. 
%%defined by $\omega\left \begin{pmatrix}
%%    x_1\\ x_2
%%\end{pmatrix}, \begin{p}
%%   y_1\\ y_2
%%\end{pmatrix})  = (x_1, x_2) \begin{pmatrix}
%%    0 & 1\\-1 & 0\end{pmatrix} \begin{pmatrix}
%%   y_1\\ y_2
%%\end{pmatrix}
%\end{enumerate}

Identifying linear endomorphism of $\mathbb{R}^2$ with two-by-two matrices, we get 
\[\mathrm{Sp}_2(\mathbb{R}) = \{ g = \begin{pmatrix}    a & b\\c & d \end{pmatrix}|\, \det{g} = ad-bc = 1\}  = \mathrm{SL}_2(\mathbb{R})\]

%Thus, let's consider $\mathbb{R}^2$. We write $\mathbb{R}$-linear endomorphisms of $\mathbb{R}^2$  as two-by-two matrices, where $ f = \begin{pmatrix}
%    a & b\\c & d\end{pmatrix}$ acts as $f((x_1, x_2)) = \begin{pmatrix}
%    a & b\\c & d\end{pmatrix}begin{pmatrix}
%    x_1\\ x_2
%\end{pmatrix} = (ax_1 + b x_2. cx_1 + dx_2)$. 
%
% and endowed with the non-degenerate skew symmetric bilinear form 
%\[
%\omega:\mathbb{R}^2 \times \mathbb{R}^2 \rightarrow \mathbb{R}, 
%\]
%defined by $\omega\left \begin{pmatrix}
%    x_1\\ x_2
%\end{pmatrix}, \begin{p}
%   y_1\\ y_2
%\end{pmatrix})  = (x_1, x_2) \begin{pmatrix}
%    0 & 1\\-1 & 0\end{pmatrix} \begin{pmatrix}
%   y_1\\ y_2
%\end{pmatrix}

\subsection{The symmetric space} 
The group $\mathrm{Sp}_2(\mathbb{R})$ acts by isometries on the hyperbolic plane $\mathbb{H}^2$, which is the symmetric space associated to  $\mathrm{Sp}_2(\mathbb{R})$. There are several different models of the hyperbolic plane, for example the Poincar\'e disk model or the upper half plane model. 

The most natural model for the symmetric space of $\mathrm{Sp}_2(\mathbb{R})$ is the space $\mathcal{J}$ of $\omega$-compatible almost complex structures on $\mathbb{R}^2$. An almost complex structure on $\mathbb{R}^2$ is a linear endomorphism, whose square is $-\mathrm{id}$. An almost complex structure $J$ is $\omega$-compatible if 
\[\omega (\cdot, J \cdot): \mathbb{R}^2 \times \mathbb{R}^2 \rightarrow \mathbb{R}\]
is a positive definite symmetric bilinear form. 

If we consider  $\mathbb{C}^2 = \mathbb{R}^2 \otimes \mathbb{C}$ the complexification of $\mathbb{R}^2$, any almost complex structure $J$ can be extended $\mathbb{C}$-linearly to an endomorphism $J_\mathbb{C}$ of 
 $\mathbb{C}^2$. This endomorphism now has two eigenspaces, $V_J^+$ and $V_J^-$ with eigenvalue $+i$ and $-i$ respectively. The symplectic form $\omega$ extends as well $\mathbb{C}$-linearly to a symplectic form $\omega_\mathbb{C}$, and gives rise to a non-degenerate Hermitian form $h(\cdot, \cdot) = i \omega_\mathbb{C} (\overline{\cdot}, \cdot)$, where $\overline{\cdot}$ is the complex conjugation. The Hermitian form $h$ is of signature $(1,1)$, and $h$ restricted to $V_J^+$  is positive definite, $h$ restricted to $V_J^-$ is negative definite. This provides two embeddings $i^{\pm}: \mathcal{J} \rightarrow \mathbb{C}P^1$, with $\overline{i^+(J)} = i^-(J)$ and $i^+(\mathcal{J}) \cup i^-(\mathcal{J})  = \mathbb{C}P^1 \backslash \mathbb{R}P^1$.

 The Poincar\'e model and the upper half plane model arise now as the image of $i^+(\mathcal{J})$ with respect to two different affine charts for $\mathbb{C}P^1$. For the Poincar\'e model consider the affine chart with respect to $V_{J_0}^-$, where $J_0$ is the standard complex structure $J_0 (x_1, x_2) = (-x_2, x_1)$. For the upper half plane model consider the affine chart with respect to the complex line spanned by the vector $e_2$. 
 
\subsection{Trace Identities}
There are two invariant functions of two by two matrices the trace ($\trace$) and the determinant ($\det$). By the classical Cayley Hamilton theorem any two by two matrix satisfies the relation 
\[ M^2 - \trace(M)M + \det(M)\mathrm{id} = 0.\] 
%Importantly, both of these functions are conjugation invariant.

%The ring of $2\times 2$ real matrices, $\Mat_2(\R)$ is one of the most ubiquitous groups in mathematics. The two functions, trace ($\trace$) and determinant ($\det$), defined from this ring to $\R$ are the most fundamental tools for analyzing matrices. This is showcased by the classical Cayley Hamilton theorem, for any $X\in \Mat_2(\R)$ 
%$$ X^2 - \trace(X)X + \det(X)\id = 0.$$ 
%Importantly, both of these functions are conjugation invariant.

Restricting to $M \in\mathrm{Sp}_2(\mathbb{R})$, where the determinant is $1$, taking the trace of the Cayley Hamilton relation, one obtains the trace identity 
\[ \trace(M^2) = \trace(M)^2 -2 .\] This can be generalized and leads to the well known $\mathrm{SL}_2$-trace identity, going back to Fricke and Klein's work \cite{Fricke-Klein-traces}: 

For any $M, N \in \mathrm{Sp}_2(\mathbb{R})$ we have 
\[ \trace(MN)+ \trace(MN^{-1}) = \trace(M)\trace(N).\]
and for the commutator:  
\begin{equation}\label{eqn:TraceCommutator}
   \trace(MNM^{-1}N^{-1}) = \trace(M)^2+\trace(N)^2+\trace(MN)^2 - \trace(M)\trace(N)\trace(MN) - 2 
\end{equation}

Trace identities have played a major role in the study of discrete subgroups of $\mathrm{SL}_2(\mathbb{R})$ and the study of representation varieties. Since the trace is a conjugation invariant function, when $S$ is a cylinder, $\pi_1(S) = \mathbb{Z}$, any representation $\pi_1(S) \to \mathrm{SL}_2(\mathbb{R})$ is essentially determined by the trace of a generator. For general surfaces, the traces of generators and some products of generators of the fundamental group give explicit coordinates for the character variety \cite{Procesi-traces}.

%In the case of a punctured torus, $\pi_1(S) = \mathbb{Z}*\mathbb{Z}$, these representations are determined by the traces of the images of the two generators and their product.  The trace of their commutator is then determined by the trace identity for the commutator.

\subsection{Markov numbers}
The trace identities allow to connect representations $\pi_1(S_{1,1}) \to \mathrm{SL}_2(\mathbb{R})$ that arise from hyperbolic structures on the once punctured torus $S_{1,1}$ to Markov numbers.

The fundamental group $\pi_1(S_{1,1})$ is a free group on two generators $A$ and $B$. 
Any representation $\pi_1(S_{1,1}) \to \mathrm{SL}_2(\mathbb{R})$ (up to conjugation) is completely determined by three coordinates: $a = \trace(A)$, $b = \trace(B)$, $c=\trace{AB}$. To make the connection with Markov numbers, one considers only representations which arise as holonomy representations of complete hyperbolic structures on $S_{1,1}$ of finite volume. For such representations the trace of the commutator $[A,B]$ is $-2$.  

Therefore the trace identitiy becomes 
\begin{equation}\label{eqn:TraceMarkov}
    a^2 + b^2+c^2 - abc = 0. 
\end{equation}
This relates the study of hyperbolic structures on a punctured torus to the classical theory of Markov triples \cite{Markoff-Markoff}, which are of interest in number theory. Markov triples are integer solutions to 
\begin{equation}
    a^2+b^2+c^2- 3abc =0
\end{equation}
Note that these two equations are related, if $(3a,3b,3c)$ is a solution to the first equation, then $(a,b,c)$ is a solution to the second equation. 

An initially surprising fact is that given a Markov triple $(a,b,c)$ the triple $(\frac{b^2+c^2}{a},b,c)$ is also a Markov triple. As proved by Markov, every Markov triple can be obtained from the ``trivial'' solution $(1,1,1)$ by repeatedly applying this rule at each entry.

Markov numbers are well studied in number theory \cite{Aigner-Markov100Years,GS-IntegralPointsMarkov}. 

\subsection{Cluster algebras}
%The structure of Markov numbers can be obtained from the cluster structure on a once punctured torus. 
A cluster algebra is an algebra with a rich additional combinatorial structure. Informally, the algebra can be generated from an initial `seed' via a process called `mutation' which produces new seeds from old seeds. The family of seeds is called cluster structure of the algebra.
Let us shortly recall the cluster structure on the once punctured torus $S_{1,1}$, which will give us a link to Markov numbers.

  \begin{equation}
        x_\gamma x_{\gamma'} = x_\alpha x_\eta + x_\beta x_\delta
    \end{equation}

\begin{figure}[hb]
    \centering
    \begin{tikzpicture}
        \node[] at (4,4) (2) {2};
        \node[] at (0,0) (4) {4};
        \node[] at (0,4) (1) {1};
        \node[] at (4,0) (3) {3};

        %Square
        \path[] (4) edge [-] node[left] {$\delta$}  (1);
        \path[] (1) edge [-] node[below] {$\gamma$} (3);
        \path[] (1) edge [-] node[above] {$\alpha$} (2);
        \path[] (2) edge [-] node[right] {$\beta$} (3);
        \path[] (3) edge [-] node[below] {$\eta$} (4); 
    \end{tikzpicture}\hspace{3pc}
    \begin{tikzpicture}
        \node[] at (4,4) (2) {2};
        \node[] at (0,0) (4) {4};
        \node[] at (0,4) (1) {1};
        \node[] at (4,0) (3) {3};

        %Square
        \path[] (4) edge [-] node[left] {$\delta$}  (1);
        \path[] (2) edge [-] node[below] {$\gamma'$} (4);
        \path[] (1) edge [-] node[above] {$\alpha$} (2);
        \path[] (2) edge [-] node[right] {$\beta$} (3);
        \path[] (3) edge [-] node[below] {$\eta$} (4); 
    \end{tikzpicture}
    \caption{Mutation at the arc $\gamma$.}
    \label{fig:QuadrilateralFlip}
\end{figure}
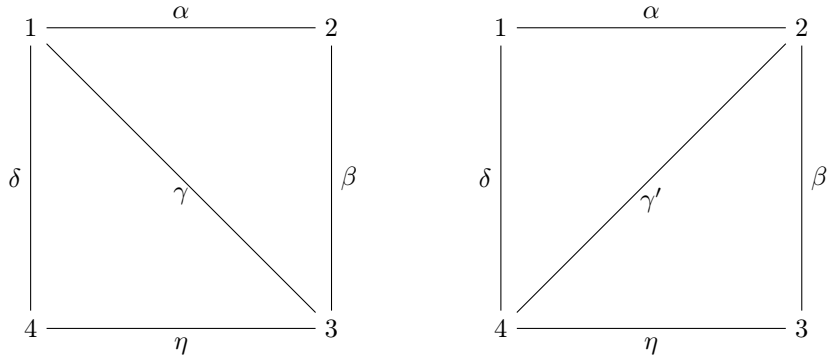

   A seed for the cluster structure of the surface $S_{1,1}$ consists of an ideal triangulation $T$ of $S_{1,1}$ along with a choice of cluster variable $x_\gamma$ for each arc $\gamma$ of $T$.   
    Each seed $(T,\vec{x})$ can be mutated at an arc $\gamma$ to produce a new seed $\mu_\gamma(T,\vec{x}) = (T',\vec{x'})$. The triangulation $T'$ is obtained by removing $\gamma$ and replacing it with the arc $\gamma'$ corresponding to the opposite quadrilateral containing $\gamma$ (\ref{fig:QuadrilateralFlip}). We name the edges of the bounding quadrilateral cyclically by $\alpha,\beta,\eta,\delta$. The cluster variables in $\vec{x'}$ are the same as the variables of $\vec{x}$ except $x_\gamma$ is replaced by $x_{\gamma'}$ which satisfies the {exchange relation}:
      The cluster algebra associated to the surface $S_{1,1}$ with initial seed $(T,x_1, \dots x_{|T|})$ is the subalgebra of the ring of rational functions in the variables of the initial seed, $\mathbb{Q}(x_1,\dots, x_{|T|})$, generated by all cluster variables obtained by repeated mutation from the initial seed.

The cluster algebra associated to a triangulated surface can be geometrically realized as the coordinate ring of decorated Teichm\"uller space. 
The decorated Teichm\"uller space of $S_{1,1}$ is the moduli space of marked complete finite volume hyperbolic metrics on $S_{1,1}$  along with a decoration consisting of a choice of horocycle around each cusp. 
The cluster coordinate $a_\gamma$ associated to an arc $\gamma$ is the $\lambda-$length of the geodesic arc in the homotopy class of $\gamma$ for the given hyperbolic metric. (See \cite{FT-LambdaLengths} for full details). Importantly, \cite{FG-DualTeichmuller} find that Penner's decorated Teichmüller space is parameterized by the positive $\mathbb{R}$-points of this cluster algebra.

\subsection{Cluster algebra and Markov triples} 
Considering the cluster structure on the once-punctured torus we can now recover again the Markov equation and Markov triples. 

\begin{figure}[hb]
    \centering
    \begin{tikzpicture}
    \node[] at (4,4) (2) {3};
\node[] at (0,0) (1) {1};
\node[] at (0,4) (3) {2};
\node[] at (4,0) (4) {4};
\path[] (2) edge [<-] node[below] {$A$}  (1);
\path[] (3) edge [<-] node[below] {$C$} (2);
\path[] (1) edge [<-] node[right] {$B$} (3);
\path[] (2) edge [->] node[left] {$B$} (4);
\path[] (4) edge [->] node[above] {$C$} (1);
\end{tikzpicture}
   \hspace{2pc}
%    \begin{tikzpicture}[]
% \node[circle,draw] at (0, 0) (2) [left
%  ] {$A$} ;
% \node[circle,draw] at (2, 3.464) (1) [] {$B$} ;
% \node[circle,draw] at (4,0) (3) [] {$C$} ;
% \path[{Square[normal,open,length=3mm]}-{Latex[transpose,open,length=3mm]}, shorten >=1.5mm,shorten <=1.5mm] (2.90) edge [] node {} (1.-150);
% \path[{Square[normal,length=3mm]}-{Latex[transpose,length=3mm]}, shorten >=1.5mm,shorten <=1.5mm] (2.30) edge [] node {} (1.-90);
% \path[{Square[normal,open,length=3mm]}-{Latex[transpose,open,length=3mm]}, shorten >=1.5mm,shorten <=1.5mm] (1.-30) edge [] node {} (3.90);
% \path[{Square[normal,length=3mm]}-{Latex[transpose,length=3mm]}, shorten >=1.5mm,shorten <=1.5mm] (1.-90) edge [] node {} (3.150);
% \path[{Square[normal,open,length=3mm]}-{Latex[transpose,open,length=3mm]}, shorten >=1.5mm,shorten <=1.5mm] (3.210) edge [] node {} (2.-30);
% \path[{Square[normal,length=3mm]}-{Latex[transpose,length=3mm]}, shorten >=1.5mm,shorten <=1.5mm] (3.150) edge [] node {} (2.30);
% \end{tikzpicture}
     \caption{The Markov quiver and a seed associated to a triangulation of a punctured torus}
     \label{fig:torusQuiver}
\end{figure}
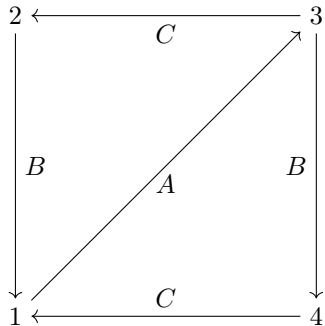

Each Markov triple corresponds to a seed of this cluster algebra (\ref{fig:torusQuiver}). The cluster mutation rule for this seed gives that the mutated coordinate $a'$  satisfies the relation $a a' = b^2 + c^2$. This is exactly the mutation of Markov triples described above. There is a mutation invariant function on this cluster algebra $f(a,b,c) = 2\frac{a^2+b^2+c^2}{abc}$. When starting with the triple $(1,1,1)$ we see that $f(1,1,1) = 6$. Since this function is fixed by mutation, every other triple obtained by mutation also satisfies $f(a,b,c) = 6$, which is equivalent to the Markov equation. The fact that all triples obtained by mutation are positive integers is a consequence of the positivity of the cluster exchange relation.

 There is also a geometric interpretation of the mutation invariant function as the sum of ``angles'' 
 \[t_i^{jk} = \frac{a_{jk}}{a_{ij}a_{ki}}\] in the triangulation. Each of these angles is the length of the horocycle at vertex $i$ inside the triangle whose edges are labeled with the cluster variables $a_{ij},a_{jk},a_{ki}$. On the punctured torus, the total length of the horocycle is then given by the Markov function:
\begin{equation}
    f(a,b,c) = 2\left(\frac{a}{bc}+\frac{b}{ac}+\frac{c}{ab}\right) = 2\frac{a^2+b^2+c^2}{abc}
\end{equation}
The length of the horocycle is clearly independent of the chosen triangulation. Furthermore the mutation of the Markov quiver leaves the adjacency of the triangles fixed. Thus the exact form of the function $f$ is invariant under mutation.

\begin{figure}
\includegraphics[scale=0.7]{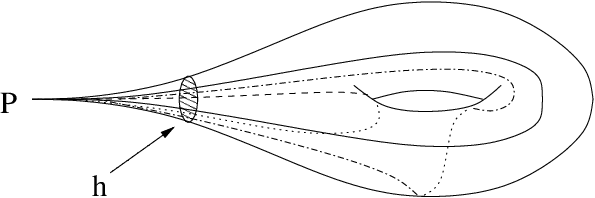}
 \caption{An ideal triangulation of a punctured torus with a horocycle. Picture taken from \cite{Teschner}.}
     \label{fig:complex}

\end{figure}

\section{Symplectic groups over non-commutative rings}
We saw above that the $\Theta$-Weyl group plays a key role in understanding the positive structures with respect to $\Theta$. In fact, a lot of the structure we see in \cite{GW_pos1, GW_posneu} suggests that if $G$ admits a positive structure with respect to $\Theta$ it might heuristically be seen as a group if type $W(\Theta)$ over a (partly) non-commutative ring. 

For Hermitian Lie groups of tube type the $\Theta$-Weyl group is of type $A_1$, and thus we would like to see a Hermitian Lie group of tube type as a group of type $\mathrm{Sp}_2$ over a non-commutative ring. This leads to the definition of groups $\mathrm{Sp}_2(A,\sigma)$ over non-commutative rings $A$ with involution $\sigma$ in \cite{ABRRW} and \cite{Rogozinnikov_20}, and in fact many of the things we described in Section~\ref{sec:SL2} can be generalized in this non-commutative setting. We first describe the construction of the groups and their symmetric space, before we go to discuss the relation with non-commutative cluster algebra, introduce non-commutative generalizations of Markov numbers, and construct  $\mathrm{Sp}_2(A,\sigma)$ quantum groups. 

We start with an involutive ring $(\A, \sigma)$, where $\A$ is a (usually) non-commutative ring with unit and where every left-inverse is also a right-inverse and $\sigma \colon \A \rightarrow \A$ an anti-involution, i.e. $ \sigma(AB) = \sigma(B) \sigma(A)$, and $\sigma(\sigma(A)) = A$ for all $A,B \in \A$. We denote by $\A^\sigma$ the set of fixed points of $\sigma$. 
A key example of an involutive ring is $\mathrm{Mat}_n(k)$ the ring of $n\times n$ matrices over a field $k$ with transpose as the anti-involution. Here the set of fixed points is the set of symmetric matrices. 
When $k = \mathbb{K} = \mathbb{C}, \mathbb{H}$, and $\sigma$ is the conjugate transpose, we have  $\A^\sigma = \Herm(n,\KK)$.

We denote by $\mathrm{Mat}_2(\A,\sigma)$ the ring of $2\times 2$ matrices over $\A$. 

We consider $\A^2 = \A\times \A$ as an $\A$-module. Given $$\matSPgreek{\Omega}_1 = \begin{bmatrix}0 & 1 \\ -1 & 0\end{bmatrix} \in \Mat_2(\A,\sigma),$$
we can define an $\A$-valued skew symmetric right-bilinear form $\omega\colon \A^2\times \A^2 \rightarrow \A$ by 
\[\omega(\vec{v},\vec{w}) = \sigma(\vec{v})^T\matSPgreek{\Omega}_1 \vec{w}.\] This form is right bilinear in the sense that 
\[\omega(\vec{v}B,\vec{w}A) = \sigma(B)\omega(\vec{v},\vec{w})A\] for $A,B \in \A$.

The symplectic group is given by the elements of $\Mat_2(\A)$ which preserve this form. 

\begin{definition}\label{def:symplectic}
The symplectic group $\SPA$ is the set of right $\A$-linear endomorphisms of $\A^2$ preserving $\omega$. 
\[\SPA = \{f\in \End(A^2) \st \forall \vec{v},\vec{w} \in \A^2: \,  \omega(f(\vec{v}),f(\vec{w})) = \omega(\vec{v},\vec{w})\}\]
In matrix form this condition is equivalent to:
\[\SPA = \{\matSP{M} \in \Mat_2(\A) \st \sigmaT{\matSP{M}}^T\matSPgreek{\Omega}_1 \matSP{M} = \matSPgreek{\Omega}_1\}. \] 
\end{definition}

When $\A =  \mathbb{R}$ and $\sigma$ is the trivial involution, then $\SPA \simeq \mathrm{Sp}_2(\mathbb{R})$. 
More generally 
when  $\A = \Mat_n(k)$ is the ring of $n\times n$ matrices over a field $k$, and $\sigma$ the matrix transpose, we have $\SPA \simeq \SP(2n,k)$. 
In fact, the classical Hermitian Lie groups of tube type can be realized essentially as $\SPA$ for appropriate choice of rings, 
see \cite{ABRRW, Rogozinnikov_20}, for the indefinite orthogonal groups $\mathrm{SO}(2,n)$ a small twist is necessary. 
When $\A =  \ \Mat_n(\CC)$ and $\sigma$ is the conjugate transpose, then $\SPA \simeq \mathrm{U}_{(n,n)}$, and if $\A =  \ \Mat_n(\HH)$ and $\sigma$ is the conjugate transpose, 
then $\SPA \simeq \mathrm{SO}^*_{(4n)}$. 
To realize the family of Hermitian Lie groups of tube type which are isogeneous to $\mathrm{SO}(2,q)$, one has to adapt the construction a bit. In this case $\A$ is the Clifford algebra associated to $\mathbb{R}^{1,q-1}$, 
and we have to consider a subgroup of $\SPA$, where matrix coefficients take values in the Clifford group.  
For the exceptional group one would like to take $\A =  \ \Mat_3(\OO)$, but this is not a ring. Thus we cannot describe the group as $\SPA$, but most of the structural aspects of the group and the Lie algebra still survive and can be described nicely in terms of Euclidean Jordan algebras. 

Note that the choice of anti-involution is crucial and dramatically affects the structure of $\SPA$ even if the underlying algebra is commutative. 
 Consider $\A = \mathbb{C}$, the complex numbers.
 In the case with trivial involution the group $\SPA$ is simply $\SL_2(\mathbb{C})$ but in the case that the $\sigma$ is complex conjugation the group $\SPA$ can be computed to be $\SL_2(\R)\times \text{U}(1)$

Above we described the symmetric space associated to $\mathrm{Sp}_2(\mathbb{R})$ as the space of compatible almost complex structures on $\mathbb{R}^2$, and derived several other models for the hyperbolic plane from that. Analogous constructions can be made for the groups  $\SPA$, starting with complex structure on $\A$-modules, where a complex structure on an right $\A$-module $V$ is an $\A$-linear map $J\colon V\to V$ such that $J^2=-\mathrm{id}$. The space of $\omega$-compatible almost complex structures on $\A^{2}$ then gives a model of the symmetric space associated to $\SPA$. From this model, one derives several other models, a projective model, a bounded (or precompact) model and an upper half space model. Thus we realize the symmetric space as a hyperbolic plane over a non-commutative ring with involution. For details we refer the reader to \cite{ABRRW}. 

When specifying to the rings so that $\SPA$ is isomorphic to a classical Lie groups $G$, this general construction provides nice models for the symmetric spaces of $G$. For classical Hermitian Lie groups of tube type the projective model gives the Borel embedding, the precompact model gives the bounded symmetric domain model and the upper half space model gives the Siegel upper half space model. 
But even in these cases, we obtain some new information. For example, the realization of the symmetric spaces as $\omega$-compatible almost complex structures is usually only described for $\mathrm{Sp}(2n,\mathbb{R})$, but not for $\mathrm{SU}(n,n)$ or $\SO^*(4n)$. Also for the complexifications of Hemitian Lie groups of tube type, we get new models of the symmetric spaces, that to our knowledge have not been described so far, see \cite{ABRRW} for the details.

\subsection{Trace Identities}
In \cite{GKW} we derived generalizations of the $\mathrm{Sp}_2(\mathbb{R})$-trace identities for the groups $\SPA$. 

For this first we recall that      $\begin{bmatrix}
        A &B\\C &D
    \end{bmatrix}$ is  $\SPA$ if it equations:
   \[
       \sigmaT(A)C-\sigmaT(C)A=0 \hspace{2pc} \sigmaT(B)D-\sigmaT(D)B=0 \hspace{2pc} \sigmaT(A)D-\sigmaT(C)B = 1= \sigmaT(D)A-\sigmaT(B)C
\]

we extend $\sigma$ to a map  
$\sigma \colon \Mat_2(\A,\sigma) \rightarrow \Mat_2(\A,\sigma)$, by setting 
$$\sigma\left(\begin{bmatrix}A & B \\ C & D\end{bmatrix}\right) = \begin{bmatrix}\sigma(A) & \sigma(B) \\ \sigma(C) & \sigma(D)\end{bmatrix}.$$ 
Note that the map $\sigma$ does not respect the ring structure on $\Mat_2(\A,\sigma)$ in any way. 
In particular, for $\matSP{M},\matSP{N} \in \Mat_2(\A,\sigma) $, $\sigma(\matSP{M}\matSP{N})$ is not equal to $\sigma(\matSP{N})\sigma(\matSP{M})$.
The involution interacts however with the matrix transpose $\matSP{M}^T$  \[\begin{bmatrix}A & B \\ C & D\end{bmatrix}^T = \begin{bmatrix}A & C \\ B & D \end{bmatrix}.\]
Note that when $\A$ is non-commutative, the matrix transpose also doesn't respect the multiplicative structure on $\Mat_2(\A,\sigma)$, but the composite transform $\matSP{M} \mapsto \sigmaT{\matSP{M}}^T$ is an anti-homomorphism.\\

We introduce $\A$-valued functions on $\Mat_2(\A,\sigma)$, which we call trace and determinant. When $\A$ is commutative and $\sigma =\mathrm{id}$ these are precisely the trace and determinants of $2 \times 2$ matrices. However, when $\sigma$ is not the identity and $\A$ is potentially non-commutative they have quite different properties.

The left determinant $\det_l \colon \Mat_2(\A,\sigma) \rightarrow \A$ is given by
    \[\begin{bmatrix}A& B\\ C & D\end{bmatrix} \mapsto \sigma(A)D - \sigma(C)B\]
    The {right determinant} $\det_r \colon \Mat_2(\A,\sigma) \rightarrow \A$ is given by 
    \[\begin{bmatrix}A& B\\ C & D\end{bmatrix} \mapsto A \sigmaT{D} - B\sigmaT{C}\]
Note that $\det_r(\matSP{M}) = \det_l(\sigmaT{\matSP{M}}^T)$. 
%However $\det_r(\sigma(\matSP{M}))$ has no relation to $\det_r(\matSP{M})$ or $\det_l(\matSP{M})$ and in fact $\det_r(\matSP{M})$ and $\det_l(\matSP{M})$ can both be invertible while $\det_r(\sigma(\matSP{M}))$ is not.
%\begin{definition}
    The trace $\trace\colon \Mat_2(\A,\sigma) \rightarrow\A$ is given by 
    \[ \begin{bmatrix}A& B\\ C & D\end{bmatrix} \mapsto A+D.\] 
    
 The trace function satisfies some interesting properties:   
 For all $\matSP{M},\Mat{N} \in \Mat_2(\A,\sigma)$ we have 
\[\trace(\sigma(\matSP{M}\matSP{N})) = \trace(\sigma(\matSP{N})\cdot \sigma(\matSP{M}))\]

    Note that in general $\sigma(\matSP{M}\matSP{N})$ is not equal to $\sigma(\matSP{N})\sigma(\matSP{M})$. If the involution $\sigma$ is trivial and thus $\A$ is commutative, the above identity reduces to the classical commutativity property of the trace $\trace(\matSP{M}\matSP{N}) = \trace(\matSP{N}\matSP{M})$, which further implies the conjugation invariance of the trace function  $\trace(\matSP{N}\matSP{M}\matSP{N}^{-1}) = \trace(\matSP{M})$. In general the trace function is however not invariant. 
%
%  The adjoint matrix of  $\matSP{M} \in \Mat_2(\A,\sigma)$ is given by \[ \Adj(\matSP{M}) =\begin{bmatrix}
%    \sigma(D) & -\sigma(B) \\
%    -\sigma(C) & \sigma(A)
%\end{bmatrix}.\]

If $\matSP{M} \in \SP_2(\A,\sigma)$ then $\matSP{M}^{-1} =\begin{bmatrix}
    \sigma(D) & -\sigma(B) \\
    -\sigma(C) & \sigma(A)
\end{bmatrix}$, and 
  $\trace(\matSP{M}^{-1}) = \sigma(\trace(\matSP{M}))$.  This is familiar from $\Sp_2(\mathbb{R})$.    
    
 In \cite{GKW} several trace identities are established 
 \begin{enumerate}
 \item {\bf Classical trace identity}:  For $\matSP{M},\matSP{N}\in \Mat_2(\A)$ we have 
    \[\trace(\Adj(\matSP{M})\matSP{N}) + \trace(\sigma(\matSP{M})\matSP{N})-\sigma(\trace(\matSP{M}))\trace(\matSP{N}) =0,\] 
    where the adjoint matrix of  $\matSP{M} \in \Mat_2(\A,\sigma)$ is given by \[ \Adj(\matSP{M}) =\begin{bmatrix}
    \sigma(D) & -\sigma(B) \\
    -\sigma(C) & \sigma(A)
\end{bmatrix}.\]
 \item {\bf Cayley Hamilton Theorem}: An element $\matSP{M} \in \Mat_2(\A,\sigma)$ lies in $\SP_2(\A,\sigma)^l$ if and only if it satisfies 
\[
    \sigma(\matSP{M})\matSP{M} - \sigma(\trace(\matSP{M}))\matSP{M} + \mathrm{id} = 0.\]
 \end{enumerate}

These trace identities are established in \cite{GKW} by explicit computations, but we expect that there is a more natural underlying structure, which relies on considering a non-commutative algebra $\A$ with involution $\sigma$ and modules over $\A$. 

\section{Relation to Manin matrices}
There is an interesting connection to Manin matrices, introduced by Manin in \cite{Manin-QuantumGroupsRemarks,Manin-QuantumGroupsNoncomGeometry}. They define some kind of $\mathrm{SL}_2$ over non-commutative rings and play a role in quantizations of $\mathrm{SL}_2$

     A matrix $\matSP{M} = \begin{bmatrix}
        A &B\\C &D
    \end{bmatrix} \in \Mat_2(R)$ is {Manin} if it satisfies the Manin commutativity equations:
   \[
       AC-CA=0 \hspace{2pc} BD-DB=0 \hspace{2pc} AD-CB = DA-BC
\]

Note that in the first two of Manin's commutativity equations, if we replace products $xy$ with $\sigmaT{x}y$ we get precisely first two condition for a matrix to be in $\SPA$. If we add the additional constraint to a Manin matrix that $AD-CB = 1$,  then the third equation matches up as well. 
%Since $\sigmaT{1} = 1$,  applying $\sigma$ to the third symplectic equation yields the corresponding equation $\sigmaT{D}A-\sigmaT{B}C = 1 = \sigmaT{A}D-\sigmaT{C}B$.

Manin matrices can be characterized by their inverse: 
    A matrix $\begin{bmatrix}
        A & B\\ C & D
    \end{bmatrix}$ is Manin if and only if its left inverse is given by $\frac{1}{AD-CB}\begin{bmatrix}D & -B\\ -C & A\end{bmatrix}$, 
    and also by satisfying the Cayley-Hamilton theorem. A matrix $M= \begin{bmatrix}
        A & B\\ C & D
    \end{bmatrix}$  is Manin if and only if it satisfies 
\[
   M^2 - \trace(M)M + \mathrm{id} = 0\]

If a matrix and its transpose are both Manin, then the entries are all forced to pairwise commute. In our situation this condition simply means that $\matSP{M}\in \SPA$. This makes $\SPA$ a more natural generalization of  $\mathrm{SL}_2$ over non-commutative rings than Manin matrices. 
For more similarities, we refer the reader to \cite{GKW}.

\section{Quantum groups}
The symplectic groups  $\SPA$  can be quantized, and this quantization parallels the quantization of $\SL_2$  \cite{Manin-QuantumGroupsRemarks,Manin-QuantumGroupsNoncomGeometry,Kassel-QuantumGroups}. 

Informally, the quantum $\SL_2$ is a deformation in a parameter $q$ of the algebra of matrix functions on $\SL_2$. 
Formally we fix a field $k$. Let $k_q = k(q)$ be the extension of $k$ by an element $q$. Often $q$ is taken to be an element of $k$, although this introduces additional subtlety if $q$ is a root of unity. 
    The coordinate ring of the quantum plane $\A_q$ is the algebra $k_q[x_1,x_2]$ such that $x_1$ and $x_2$ q-commute, i.e. $x_1x_2 = q^{-1}x_2x_1$.
Intuitively, the quantization of $\SL_2(k)$ should act like the endomorphisms of the quantum plane. 

The algebra of quantum matrix functions $\CoordRing(\SL_2)_q$ is the ring $k_q\{a,b,c,d\}/I$ where $I$ is generated by the quantum $\SL_2$ relations:
    \begin{align}\label{eqn:QuantumSl2}
\begin{aligned}
     ac-q^{-1}ca =& 0  \hspace{3pc}& ab-q^{-1}ba =& 0\\
     bd-q^{-1}db =& 0 & cd-q^{-1}dc =& 0\\
     ad-q^{-1}cb =& da-qbc  &  bc=cb \\
     ad-q^{-1}cb =& 1
\end{aligned}
\end{align}
These equations can be obtained by adding the appropriate power of $q$ to the equations for a matrix and its transpose to be Manin.  To better compare with $\SPA$ we can rewrite the quantum $\SL_2$ relations as follows
\begin{align}\label{eqn:QuantumSL2rephrased}
    \begin{aligned}
        ac-q^{-1}ca =& 0 \hspace{3pc} & ab-q^{-1}ba =& 0\\
        bd-q^{-1}db =& 0 & cd-q^{-1}dc =& 0\\
        ad-q^{-1}cb =& 1  & ad-q^{-1}bc =&1 \\
        da-qbc =&1 & da-qcb =&1
    \end{aligned}
\end{align}

This algebra is actually a Hopf algebra. We will not address this here, but refer the reader to \cite{Kassel-QuantumGroups} and \cite{GKW} for more details. 

To quantize $\SPA$ we first add a quantum parameter $q$ to our algebra $\A$. Formally let $\A_q = \A[q^{\pm}]$ with the additional assumption that $q$ commutes with all of $\A$. We then extend $\sigma$ to $\A_q$ by setting $\sigma(q) = q^{-1}$.

%    Two elements $X,Y \in \A_q$ \keyword{left $q-\sigma-$commute} if $\sigmaT{X}Y = q^{-1}\sigmaT{Y}X  $.

 Following the quantization of $\SL_2(\R)$ we define a system of equations on the matrix entries of a $2\times 2$ matrix specifying the quantum commutativity and quantum determinant of the elements.

     The quantum symplectic equations for a matrix $\begin{bmatrix}
         A & B\\ C &D
     \end{bmatrix} \in \Mat_2(\A_q)$ are
       \begin{align}\label{eqn:QuantumSymplectic}
    \begin{aligned}
         \sigma(A)C-q^{-1}\sigma(C)A =& 0  \hspace{3pc}& A\sigma(B)-q^{-1}B\sigma(A) =& 0\\
         \sigma(B)D-q^{-1}\sigma(D)B =& 0 &  C\sigma(D)-q^{-1}D\sigma(C) =& 0\\
         \sigma(A)D-q^{-1}\sigma(C)B =& 1 &  A\sigma(D)-q^{-1}B\sigma(C) =& 1
    \end{aligned}
    \end{align}

 Note that we consider $\sigmaT{A}D-q^{-1}\sigmaT{C}D$ to be the quantum determinant.

Let us point out that these equations correspond both to the standard symplectic equations and to the rephrased quantum $\SL_2(\R)$ equations (\ref{eqn:QuantumSL2rephrased}). We stress that these equations do not (!) imply any commutativity between $B$ and $C$. 

Usually when we define a quantum deformation of an algebra, we  are constructing a non-commutative algebra out of a commutative one. Here we already start with a non-commutative algebra, so its quantization is a much less drastic change. In fact the corresponding quantization can be seen simply as a change of symplectic form. 
The equations above are equivalent to $\matSP{M}$ having a two sided inverse and preserving the $\A_q$ valued skew symmetric right-bilinear form $\omega_q : \A_q^2 \times \A_q^2 \rightarrow \A_q$ by
\[\omega_q(\vec{v},\vec{w}) = \sigmaT{\vec{v}}^{T}\Omega_q \vec{w}, \] with \[\matSPgreek{\Omega}_q = \begin{bmatrix} 0& 1\\ -q^{-1} & 0\end{bmatrix} \in \Mat_2(\A_q).\]

 The corresponding group $\SPA_q$ acts on the left quantum plane, $(\A^2_\omega)_q$, which is the set of generic left $q$-commuting elements of $\A_q^2$.
    \[(\A^2_\omega)_q = \{\vec{v}=(V_1,V_2) \in \A_q^2 \mid \text{$\vec{v}$ generic and } q^{-1}\sigmaT{V_2}V_1 = \sigmaT{V_1}V_2\}\]

In fact, and contrary to the commutative situation, the group $\SPA_q$ has many elements if we assume that $q$ has a square root $q^{1/2}$. This is very different from the usual quantization of $\SL_2(\R)$ where the corresponding group is just $\left\{\begin{bmatrix}
       a & 0\\0 & a^{-1}
   \end{bmatrix} \mid a \in \R\right\}.$
  Here  $\matSP{M} = \begin{bmatrix}
        A & B\\ C &D 
    \end{bmatrix} \in \SPA$ then $\begin{bmatrix}
        A & q^{1/2}B\\q^{-1/2}C & D
    \end{bmatrix} \in \SPA_q$. 

One can further introduce the algebra of quantum symplectic matrix functions and endow it with an appropriate Hopf algebra like structure. For this construction we refer the reader to \cite{GKW}. There the reader will also find a discussion of the relation to the classical algebra of quantum matrix functions when $\A = k$ is a field and $\sigma$ is the trivial involution. 

\begin{question}
Is there an interpretation of this quantization in terms of double Poisson brackets? 
\end{question}

\section{Non-commutative cluster algebras}
In \cite{BR-noncommutativeClusters} Berenstein and Retakh introduced a non-commutative generalization of the $A_1$- cluster algebras associated to a surface.  

For the non-commutative $A_1$-type cluster algebra associated to a surface $S$, to construct a {seed} , we start again with an ideal triangulation of the surface $S$. Each possible orientation of an edge $i \rightarrow j$ of the triangulation is assigned a variable $X_{ij}$ in a non-commutative algebra. 
     There is an additional condition for each triangle cyclically labeled $i,j,k$. Berenstein and Retakh introduce a {non-commutative angle} $T_i^{jk} = X_{ji}^{-1}X_{jk}X_{ik}^{-1}$. The seed has to satisfy the triangle condition: $T_i^{jk} = T_i^{kj}$. 
The non-commutative cluster algebra is then again grown through  mutation. Each mutable edge $i \rightarrow k$ is contained in a unique quadrilateral with vertices $(i,j,k,\ell)$. The mutated seed contains the triangulation with $i\rightarrow k$ replaced with $j \rightarrow \ell$. The element assigned to $j\rightarrow \ell$ is given by the formula 
 \[
     X_{j\ell} = X_{jk}X_{ik}^{-1}X_{i\ell} + X_{ji}X_{ki}^{-1}X_{k\ell}
\]

% The non-commutative surface cluster algebra, $C(S)$ associated to an initial seed $(T,X_1, \dots X_{|T|})$ is the subalgebra of the non-commutative fraction field of the $\mathbb{Q}(\{X_i\})$ generated by all of the variables $X_{\gamma}$ found in all seeds modulo the triangle conditions. 
%%    
%   This algebra becomes a $(\mathbb{Q},\sigma)$-algebra once we define 
   We want to relate the variables associated to the same edge with different orientation, so we set 
   $\sigma(X_{ij})=X_{ji}$ and extend $\sigma$ as an anti-involution.  
    Using $\sigma$, we see that the triangle condition reads as $$T_i^{jk} = T_i^{kj} = \sigma(T_i^{jk}),$$ so that each non-commutative angle is in $C(S)^\sigma$ 
    
    The mutation rule is equivalent to the requirement that the angles are additive.
     \begin{align*}
         T_{i}^{j\ell} =~& T_{i}^{jk}+T_{i}^{k\ell}\\
         {X_{ji}^{-1}}X_{j\ell}{X_{i \ell}^{-1}} =~& {X_{ji}^{-1}}X_{jk}{X_{ik}^{-1}} + {X_{ki}^{-1}}X_{k\ell}{X_{i \ell}^{-1}}
     \end{align*}

The non-commutative cluster structure can also be used to parametrize the positive semigroups $U_\Theta^{>0}$ and $G_\Theta^{>0}$ for classical Hermitian Lie groups of tube type. This is discussed in \cite{GKW, GKNW1}. 

Recall that in the commutative $A_1$-cluster algebra associated to the surface $S$ we had a geometric realization by looking at parametrizations of the decorated Teichm\"uller space by Penner's $\lambda$-lengths. The decorated Teichm\"uller space can be realized as the positive points of the $\mathbb{R}$-points of the cluster algebra. 

There is also a geometric realization of the non-commutative cluster algebra, using positivity. 
In \cite{AGRW} the authors introduced non-commutative $\Lambda$-lengths for decorated representations $\pi_1(S) \rightarrow \mathrm{Sp}(2n,\mathbb{R})$ and show that these $\Lambda$-lengths give a geometric realization of the Berenstein-Retakh cluster algebra, such that the appropriate positive points are precisely the space of decorated maximal representations, i.e. the space of positive representations with respect to the positive structure on $\mathrm{Sp}(2n,\mathbb{R})$  with respect to $\Theta = \{\alpha_n\}$.  See also \cite{GKW, Rogozinnikov_20,Kineider_Rogozinnikov} for more general symplectic groups $\SPA$ over non-commutative rings, and \cite{Kineider_Rogozinnikov, Goncharov_Kontsevich} for related work for $\mathrm{GL}_n$ over non-commutative rings. 

These geometric realizations of the non-commutative cluster algebra of type $A_1$ provide the non-commutative analogue of the $\mathcal{A}$-cluster variety of Fock and Goncharov in the commutative setting. In \cite{AGRW} the authors also provide $\mathcal{X}$-type coordinates, but they are not well adapted to the cluster structures. Goncharov and Shen \cite{GoncharovShen} introduced $\mathcal{P}$-cluster varieties, which seems to be more appropriate for a non-commutative generalization. 
\begin{questions}
\begin{enumerate}
\item Provide appropriate non-commutative $\mathcal{P}$-cluster varieties of type $A_1$. 
\item In \cite{FockGoncharov} tropical points of cluster varieties play an important role. In the $A_1$-situation they parametrize laminations on the surface. Provide appropriate tropical non-commutative $\mathcal{A}$-cluster varieties and $\mathcal{P}$-cluster varieties and geometrically describe what they parametrize. 
\item For split real Lie groups, Fock and Goncharov \cite{FockGoncharov} formulated wide ranging duality conjectures, which have been proved in \cite{GrossHackingKeelKontsevich}. 
It is an open question if there are any dualities for Lie groups admitting a positive structure with respect to $\Theta$.
 \end{enumerate}
\end{questions}

\section{Non-commutative Markov numbers}
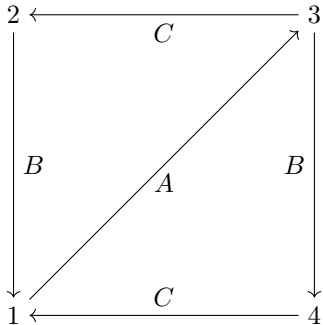
\begin{figure}[hb]
    \centering
    \begin{tikzpicture}
    \node[] at (4,4) (2) {3};
\node[] at (0,0) (1) {1};
\node[] at (0,4) (3) {2};
\node[] at (4,0) (4) {4};
\path[] (2) edge [<-] node[below] {$A$}  (1);
\path[] (3) edge [<-] node[below] {$C$} (2);
\path[] (1) edge [<-] node[right] {$B$} (3);
\path[] (2) edge [->] node[left] {$B$} (4);
\path[] (4) edge [->] node[above] {$C$} (1);
\end{tikzpicture}
   \hspace{10pc}
%    \begin{tikzpicture}[]
% \node[circle,draw] at (0, 0) (2) [left
%  ] {$A$} ;
% \node[circle,draw] at (2, 3.464) (1) [] {$B$} ;
% \node[circle,draw] at (4,0) (3) [] {$C$} ;
% \path[{Square[normal,open,length=3mm]}-{Latex[transpose,open,length=3mm]}, shorten >=1.5mm,shorten <=1.5mm] (2.90) edge [] node {} (1.-150);
% \path[{Square[normal,length=3mm]}-{Latex[transpose,length=3mm]}, shorten >=1.5mm,shorten <=1.5mm] (2.30) edge [] node {} (1.-90);
% \path[{Square[normal,open,length=3mm]}-{Latex[transpose,open,length=3mm]}, shorten >=1.5mm,shorten <=1.5mm] (1.-30) edge [] node {} (3.90);
% \path[{Square[normal,length=3mm]}-{Latex[transpose,length=3mm]}, shorten >=1.5mm,shorten <=1.5mm] (1.-90) edge [] node {} (3.150);
% \path[{Square[normal,open,length=3mm]}-{Latex[transpose,open,length=3mm]}, shorten >=1.5mm,shorten <=1.5mm] (3.210) edge [] node {} (2.-30);
% \path[{Square[normal,length=3mm]}-{Latex[transpose,length=3mm]}, shorten >=1.5mm,shorten <=1.5mm] (3.150) edge [] node {} (2.30);
% \end{tikzpicture}
     \caption{Seed associated to a triangulation of a punctured torus}
     \label{fig:torusQuiver_non}
\end{figure}

We now use the non-commutative $A_1$-cluster algebra structure to obtain non-commutative generalizations of Markov numbers. For this we replace the real numbers $\mathbb{R}$ by an involutive algebra $(\A,\sigma)$.  For details we refer the reader to \cite{GKW}.
We start with the quiver associated to the once-punctured torus, see Figure~\ref{fig:torusQuiver_non}

Now we do not have a horocycle and the length of the pieces of this horocycle to build our invariant function. We considered the lengths of horocycle arcs $t_i^{jk}$, they are now ``upgraded'' to the non-commutative angle $T_i^{jk}$, which appear in the triangle conditions of the non-commutative cluster algebra. We require these angles to be fixed by the anti-involution $\sigma$.\\
Explicitly, in $A,B,C$ are the seed variables, see Figure~\ref{fig:torusQuiver_non}, these angles are all cyclic permutations of:
\[\sigmaT{A^{-1}}B\sigmaT{C^{-1}} \in \A^\sigma \hspace{2pc} \sigmaT{A^{-1}}C\sigmaT{B^{-1}} \in \A^\sigma \]

  The {non-commutative mutation} at the node $A$ is given by the equation:
    \begin{equation}\label{eqn:MarkovMutation}
        A' = \sigma(C)A^{-1}\sigma(C) + B\sigma(A^{-1})B
    \end{equation}
The generalized Markov function, which is invariant under mutations, for a seed $(A,B,C)$ whose triangles are cyclically oriented, is  
     \begin{equation}\label{eqn:NonComMarkov}
     \begin{aligned}
         F(A,B,C) =& ~\hphantom{+}\sigma(A^{-1})B\sigma(C^{-1}) + \sigma(B^{-1})C\sigma(A^{-1})+\sigma(C^{-1})A\sigma(B^{-1})\\
         &+\sigma(C^{-1})B\sigma(A^{-1})+\sigma(A^{-1})C\sigma(B^{-1})+\sigma(B^{-1})A\sigma(C^{-1})
         \end{aligned}
     \end{equation}
The reader might check that if the involution $\sigma$ is trivial and $a,b,c$ commute, this simplifies $F$ to $F(a,b,c) = 2\frac{a^2+b^2+c^2}{abc}$ which is the standard Markov function. 

Now we  can to use the generalized Markov function to consider generalizations of Markov triples. In the classical case we get Markov triples by starting with the seed $(A,B,C) = (1,1,1)$ and by mutation we get only integer solutions. 
There are two important things we have to take into account
\begin{enumerate}
\item {\bf Integrality} 
From the cluster theory there is nothing special about the integers as a subring of $\R$. But the 
Laurent phenomenon for cluster algebras states that every number obtained by mutation can be expressed as
a Laurent polynomial in the original triple $(a, b, c)$. Thus if $(a, b, c)$ are units in any subring all values
obtained by mutation belong to that subring. 
Thus for non-commutative generalizations we are interested in starting with a seed that takes values in units of an appropriate subring of the integral  elements in $(\A,\sigma)$. Recall the integral elements of a ring $\A$ are the roots of monic polynomials with integer coefficients in $\A$.
\item{\bf Positivity} 
Another subtlety when choosing an initial triple is that we require every element of the cluster algebra is invertible. Simply choosing a triple of invertible elements isn't enough to guarantee that all elements obtained by mutation are invertible. For example in $\CC$ with the trivial involution, mutation at 1 in the triple $(1,1,i)$ results in the triple $(0,1,i)$. A sufficient solution to avoid these triples is to require that all non-commutative angles in the initial seed are ``positive''  in the appropriate sense. 

Over $\R$ it is clear what positive means, but over arbitrary involutive rings it is less clear. We define a positive structure on $(\A,\sigma)$ to be a subset $P$ of $\A^\sigma$ containing 1 and be closed under  addition, inversion, and $\sigma$-congruence ($A \mapsto \sigma(B)AB$) by invertible elements.  Here we see the connection with the positive structure on Hermitian Lie groups of tube type. 
\end{enumerate}

Thus, putting these things together, we fix a subring $R$ of integral elements of $\A$ and a positive structure $P$ on $(\A,\sigma)$. 
An admissible seed (with respect to $R$ and $P$) is then a triple $(a,b,c) \in \A^3$, where $a,b,c$ are units in $R$ and all the non-commutative angles are in $P$.

A simple way to get an admissible triple is to take a unit $B$ in $R$, then the triple $(1,B,B^{-1})$ is admissible.

With this general setup we now obtain several examples of non-commutative Markov triples. 

\subsection{Matrix Markov numbers} 
Consider $\A = \Mat(n,\R)$ with $\sigma$ being the standard matrix transpose. We choose $R = \mathrm{Mat}_n(\ZZ)$ as the subring of integral elements to focus on. The positive structure $P$ consists of the set of positive definite symmetric matrices.

There are many complicated examples of admissible triples, let us describe just the simplest one (see \cite{GKW} for more examples). 
For this let $B_t = \begin{bmatrix}
     1 & t\\ 0 & 1
 \end{bmatrix}$ and consider the seed $(1,B_t,B_t^{-1})$. Then the value of the Markov function is 
 \[\begin{bmatrix}
    6+2t^2 & 0\\
    0 & 6 + 2t^2
\end{bmatrix}\]

The matrices we obtain are matrix Markov number. 
Then $\frac{1}{2}$ times the trace of each matrix Markov number is a polynomial  with integral coefficients $p(t)$, where $p(0)$ is a classical Markov number. Thus in this case $\frac{1}{2}$ trace of each matrix Markov number can be considered as a one dimensional deformation of classical Markov numbers.

\subsection{Complex and dual  Markov numbers} 
We also obtain interesting examples for commutative rings with non-trivial involutions $\sigma$. One interesting example is to consider the complex numbers $\mathbb{C}$ with the involution $\sigma$ which is the complex conjugation. Then $\A^{\sigma} = \mathbb{R}$.  
For the ring of integral elements $R$ we take the roots of unity. The positive structure is given by the positive real numbers. 

A simple admissible seed triple is a triple $z_1,z_2,z_3\in \CC$ of three roots of unity whose product is 1. If these roots of unity are chosen to be nontrivial, the complex Markov numbers we get differ from the real Markov numbers. In this case the mutation rule splits into a mutation rule for the norm, which is just the classical mutation rule for Markov triples, and a mutation rule for the arguments. Therefore the norm of these complex Markov numbers are classical real Markov numbers. It is however interesting to understand the arguments of these complex Markov 
numbers. 
In \cite{GKW} the following conjecture is made 
\begin{conjecture}
 If the initial seed triple consists of 3 roots of unity $(z_1,z_2,z_3)$ such that $n$ is the minimal integer with $z_i^n = 1$ then every triple of $\frac{a}{n},\frac{b}{n},\frac{c}{n}$ that sums to $1$ is achieved as the argument of a complex Markov triple. Furthermore if $(z_1,z_2,z_3)$ are irrational angles then the image of the arguments of Markov triples is dense in the circle.
\end{conjecture} 

\begin{figure}
\includegraphics[scale=0.6]{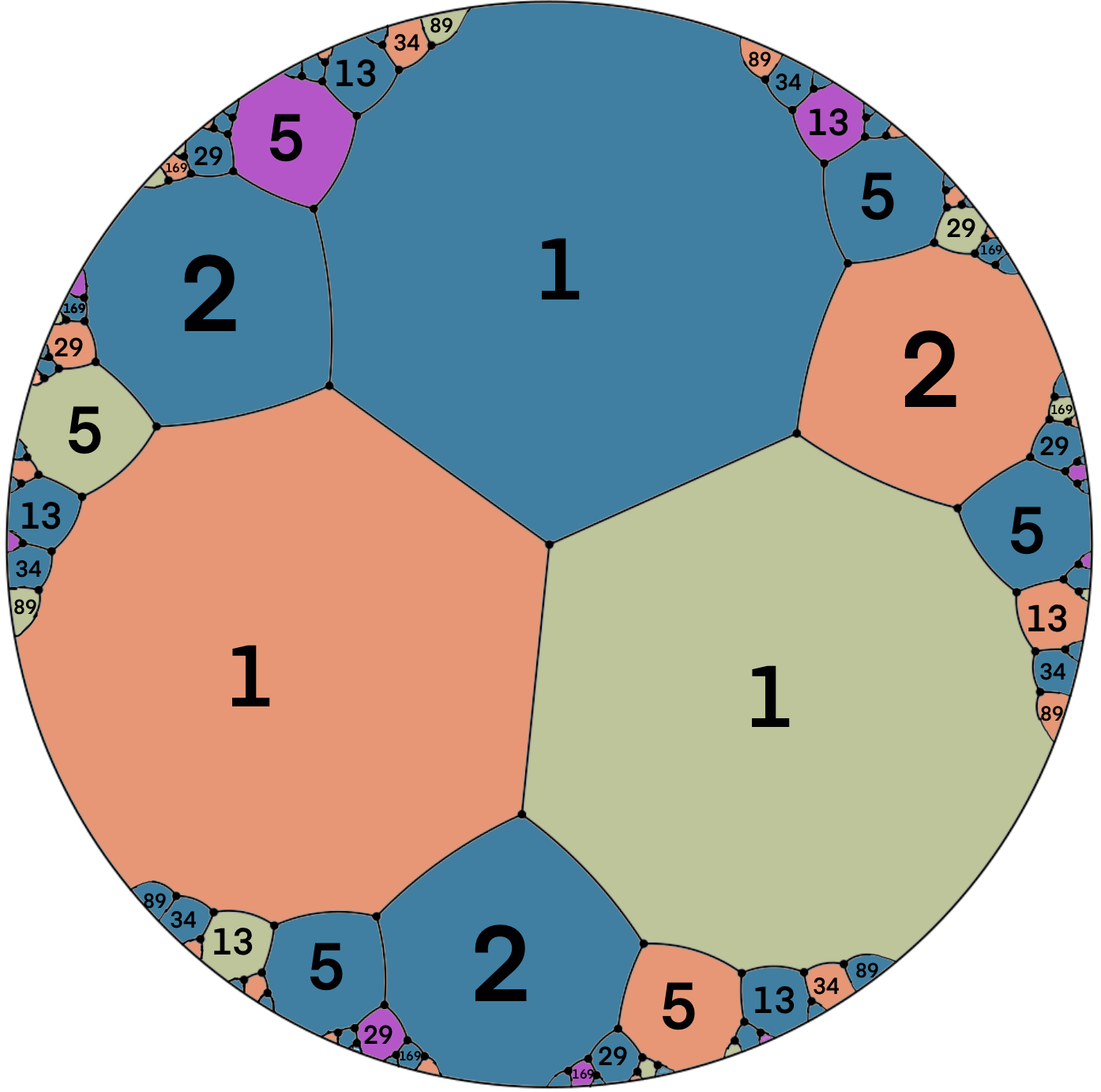}
 \caption{The Markov ``tree'' for complex Markov numbers}
     \label{fig:complex}

\end{figure}

Another interesting example to consider are the dual numbers $\{a + \alpha\epsilon | a,\alpha\in \RR, \epsilon^2 = 0\}$, and the subring $R$ of integral elements of the form $a + \alpha \epsilon$ with $a,\alpha \in \ZZ$. The units of this subring are elements of the form $1 + \alpha\epsilon$ for $\alpha \in \ZZ$. We can write dual numbers in the form $a(1+\alpha\epsilon)$ with $a \neq 0$. This is analogous to the polar form of a complex number.

If one considers the trivial involution, one obtains dual Markov numbers which have been studied also \cite{BoninOvsienko}. The more interesting case is again when $\sigma$ is the non-trivial involution 
given by  the conjugation map $a+\alpha\epsilon \mapsto a -\alpha \epsilon$. In this case the mutation rule splits again, see \cite{GKW}

%\section{Non-commutative Algebra}
\section{Beyond the $A_1$-case}
In the above sections we illustrated how the non-commutative perspective suggested by positivity leads to a fresh perspective on a variety of things. 
We mostly focused on the $A_1$-case which arises from the positive structure on Hermitian Lie groups of tube type with respect to a subset $\Theta$ of $\Delta$ consisting of one element. 
This already extends beyond just the family of Hermitian Lie groups of tube type. When considering the groups $\SPA$, the trace identities, their quantizations or the corresponding Markov triples,  the non-commutative algebra $\mathcal{A}$ does not have to be of Hermitian type  - any (non-commutative) ring $\mathcal{A}$ with ant-involution $\sigma$ works. And actually, one does not have to stop at skew-symmetric bilinear forms, but can also consider symmetric bilinear forms over $(\mathcal{A}, \sigma)$, see \cite{ABRRW}.

The phenomena we describe in fact suggest that there is a more fundamental approach to explain some of the observations by considering $\sigma$-commutative modules over $(\mathcal{A}, \sigma)$. More on this perspective will appear in forthcoming work with Dani Kaufmann. 

Going back to the families of simple Lie groups admitting a positive structure, the next family is formed by the indefinite orthogonal groups (or indefinite spin groups), whose $\Theta$-Weyl group is of type $B_p$. Here the positive structure suggests to look for a partially non-commutative structure.  
One can in fact introduce partially non-commutative cluster structures of type $B_p$. This leads to the introduction of non-commutative polygonal cluster algebras \cite{GKNW1}. 

In \cite{GKNW2}, building on work of Goncharov and Shen \cite{GoncharovShen} we introduce more general partially non-commutative cluster structures associated to more general $\Theta$-Weyl groups. These cluster structures can be used to parametrize the open semigroups $U_\Theta^{>0}$ and $G_\Theta^{>0}$, to give explicit positivity tests, and provide (partially) non-commutative cluster structure on the space of decorated positive representations. 

%------
% Insert acknowledgments and information
% regarding funding at the end of the last
% section, i.e., right before the bibliography.
%------

%\begin{ack}
%I thank my collaborators for taking this mathematical journey with me. 
%\end{ack}
%
%\begin{funding}
%This work  was partially  supported by the European Research Council under ERC-Advanced Grant 101018839, and by the Deutsche Forschungsgemeinschaft  TRR 191. The author thanks the Hector Fellow Academy for support.
%\end{funding}

%------
% Insert the bibliography.
%------

\bibliography{ECM.bib}
\bibliographystyle{alpha}

%\begin{thebibliography}{99}

%------ Example for a paper in journal:
% \bibitem{article1}
% A.~Petrunin, Parallel transportation for Alexandrov space with curvature bounded below.
% \emph{Geom. Funct. Anal.} \textbf{8} (1998), no.~1, 123--148
% \Zbl{0903.53045} \MR{1601854}

%------ Example for a book:
% \bibitem{book1}
% W.~P. Ziemer, \emph{Weakly differentiable functions}.
% Grad. Texts in Math. 120,  Springer, New York, 1989
%\Zbl{0692.46022} \MR{1014685}

%------ Example for a paper in a book:
% \bibitem{incollection1}
% J.~S. Milne, Introduction to Shimura varieties.
% In \emph{Harmonic analysis, the trace formula, and Shimura varieties},
% pp. 265--378, Clay Math. Proc. 4,
% American Mathematical Society, Providence, RI, 2005
% \Zbl{1148.14011} \MR{2192012}

%------ Example for a preprint on arXiv:
% \bibitem{preprint1}
% D.~V. Nguyen, S.~K. Chilappagari, M.~W. Marcellin, and B.~Vasic,
% LDPC codes from latin squares free of small trapping sets.
% 2010, \arxiv{1008.4177}

%------ Example for a report:
% \bibitem{report1}
% J.~Schöberl, Commuting quasi-interpolation operators.
% Technical report isc-01-10-math, Texas A\&M University, 2001,
% \url{www.isc.tamu.edu/publications-reports/tr/0110.pdf}

%------ Example for a thesis:
% \bibitem{thesis1}
% E.~Giorgi, \emph{The geometric universe}.
% Ph.D. thesis, University of Maryland, College Park, 2002

%\end{thebibliography}

\end{document}